\documentclass[a4paper,11pt,reqno]{article}
\usepackage[T1]{fontenc}
\usepackage{courier}
\usepackage[scaled=0.92]{helvet}
\usepackage{amscd,amsmath,amssymb,amsthm,amsfonts,epsfig,graphics,accents}

\usepackage[percent]{overpic} 				
\usepackage[english]{babel}
\usepackage{algpseudocode,algorithm,float}
\usepackage{amsthm,amsfonts,amssymb,amsmath,subfigure}
\usepackage{hyperref}
\hypersetup{
colorlinks=true,
citecolor=red,
linkbordercolor={1 1 1},
citebordercolor={1 1 1},
}
\usepackage{cleveref}
\usepackage{mathrsfs}
\usepackage{color}
\usepackage{bm}

\usepackage{xcolor} 

\newcommand{\mx}{\mbox}
\newcommand{\rw}{\rightarrow}
\newcommand{\de}{\displaystyle}
\newcommand{\ml}{\mathcal}

\newcommand{\pl}{\partial}

\newcommand{\x}{\times}

\newcommand{\beq}[1]{\begin{equation} \label{#1}}
\newcommand{\eeq}{\end{equation}}
\newcommand{\beqar}{\[ \begin{array}{rcl}}
\newcommand{\eeqar}{\end{array} \]}

\newcommand{\ue}[1]{\boldsymbol{#1}}
\newcommand{\uue}[1]{{\boldsymbol{#1}}}

\providecommand{\ep}{\varepsilon}
\providecommand{\ph}{\varphi}

\providecommand{\RR}{\mathbb{R}}

\providecommand{\NN}{\mathbb{N}}

\newtheorem{satz}{Theorem}[section]
\newtheorem{prop}{Proposition}[section]

\newtheorem{lem}[satz]{Lemma}

\newtheorem{defn}[satz]{Definition}
\newtheorem{rem}{Remark}[section]

\newtheorem{hyp}{Hypothesis}[section]

\usepackage{xcolor}

\usepackage{xcolor}

\providecommand{\ep}{\varepsilon}
\providecommand{\ph}{\varphi}

\providecommand{\RR}{\mathbb{R}}

\providecommand{\NN}{\mathbb{N}}

\usepackage[utf8]{inputenc}

\hypersetup{pdfauthor={Fortunati, A. et al.},%
            pdftitle={Analytical solutions...},%
            pdfsubject={},%
            pdfkeywords={}%
}

\DeclareMathOperator{\diag}{diag}

\DeclareMathOperator{\flag}{flag}
\DeclareMathOperator{\true}{true}

\makeatletter
\renewcommand*{\@fnsymbol}[1]{\ensuremath{\ifcase#1\or *\or \mathsection \or \dagger \or \else \fi}}
\makeatother

\title{\LARGE{\textbf{{Analytical solutions for ultra-fast precessional switching in inertial magnetization dynamics}}}}
\date{}
\author{%
Alessandro Fortunati\thanks{E-mail: alessandro.fortunati@unina.it}, Massimiliano d'Aquino\thanks{E-mail: mdaquino@unina.it} $\,$ and Claudio Serpico\thanks{E-mail: serpico@unina.it}   
\bigskip \\
Department of Electrical Engineering and Information Technologies, \\University of Naples Federico II, I-80125 Naples, Italy.
}

\usepackage{amsmath}
\usepackage{graphicx}
\begin{document}

\maketitle

\begin{abstract}
Here we consider the magnetization dynamics in ferromagnetic nanoparticles or films driven by external magnetic field pulses directed transverse to the initial equilibrium state. The excitation pulse drives large-angle ultra-fast magnetization dynamics that may eventually end up in the reversed equilibrium realizing successful precessional switching. We consider ultra-short pulse duration (fractions of picosecond) and large external field intensities (several Tesla) which may be relevant for the realization of faster energy-efficient memory cells. For such short time scales, we include inertial effects in the theoretical description by considering the inertial Landau-Lifshitz-Gilbert equation which requires to be treated as a system of singularly perturbed ODEs for small values of the inertia. By using suitable perturbation approach based on multiple time scales analysis, we develop an approximate closed-form solution for the switching dynamics as well as formulas for switching time and its safety tolerances to obtain the successful switching as function of physical parameters of the particle and strength of magnetic inertia. The developed analytical solution is validated by numerical simulations.

 
\end{abstract}


\section{Introduction}

The study of ultra-fast magnetization dynamics is a central theme both for the advancement of knowledge in the area of applied magnetism and for its important implications on the realization of improved magnetic and spintronic nanotechnologies\cite{Dieny2020}, such as faster magnetic recording technologies (e.g. magnetic random access memory - MRAM - cells, advanced sensors or non boolean logic devices).  

 With the ever-growing demand for faster and more energy-efficient data storage, understanding and controlling magnetization switching at picosecond and even shorter timescales becomes essential. Traditional magnetic recording methods or those which rely on thermal activation (e.g. heat-assisted magnetic recording - HAMR) are inherently limited in terms of speed and energy efficiency. In contrast, precessional switching, where magnetization is driven by strong and ultra-short external magnetic field\cite{kaka2002precessional} or spin-torques pulses\cite{bedau2010spin}, has potential to overcome these barriers by enabling ultra-fast, non-thermal switching mechanisms.

In spite of the potential advantages in terms of speed, precessional switching suffers from the strong sensitivity to the timing of the external excitation pulse, which must be very precise in order to switch off at the right instant and let magnetization relax to the target reversed orientation. 
This issue has been the focus of theoretical studies\cite{bertotti2003geometrical,daquino2004numerical,BMS2009,fortunati2024} that provided analytical formulas for switching time as function of material and geometrical parameters of the nanomagnet. 

\begin{figure}
    \centering
    \includegraphics[width=0.5\linewidth]{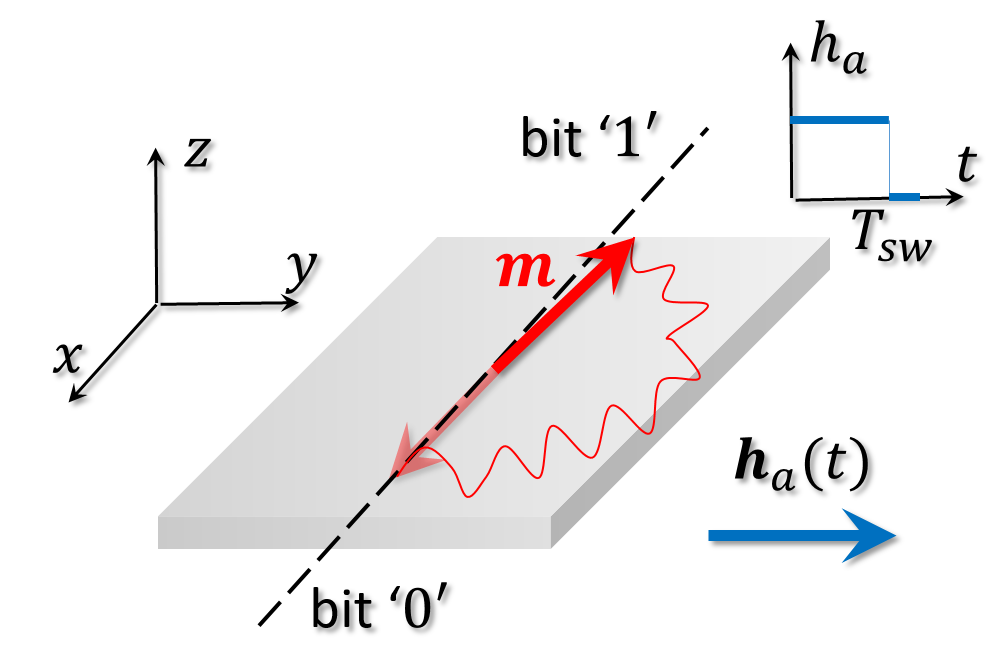}
    \caption{Sketch of ultra-fast precessional magnetization switching via ultra-short transverse external field pulses. When the duration $T_{sw}$ is below picoseconds, inertial effects set in appearing as THz nutation (thin red line) superimposed to magnetization precession.}
    \label{fig:prec switching sketch}
\end{figure}

However, when ultra-short (below ps) external pulses with large intensities (several Tesla) are used, the emergence of inertial effects on magnetization dynamics becomes relevant in that magnetization undergoes nutations at THz frequencies superimposed to classical precessional dynamics (see Fig.\ref{fig:prec switching sketch} for a sketch), as it was experimentally demonstrated in recent times\cite{Neeraj2021inertial} after the theoretical prediction one decade before\cite{Ciornei2011magnetization}. For this reason, in theoretical analysis of such processes, one has to consider the inertial version of the Landau-Lifshitz-Gilbert (LLG) equation, termed hereafter iLLG equation\cite{Ciornei2011magnetization,neeraj2022inertial}.

Despite its apparent similarity with the classical LLG, the iLLG equation has profoundly different mathematical nature owing to the presence of the second-order time derivative of magnetization that implies the definition of a higher dimensional nonlinear dynamical system. Nevertheless, such equation has remarkable qualitative properties related with conservation of quantities and suitable energy balance\cite{DAQUINO2024112874}. As a result, for instance in the general spatially-inhomogeneous magnetization case, this equation describes propagation of short wavelength spin-waves (tens of nanometer) with velocities that significantly deviate from the classical behavior due to the onset of an ultimate speed limit in the range of several thousand meters per second\cite{daquino2023micromagnetic}. 

In this paper, we study the iLLG dynamics for a macrospin particle driven by intense external field pulses that are transverse to the initial equilibrium magnetization in order to provide approximate solutions for the temporal dynamics of the magnetization (unit) vector field and give closed-form formulas for switching time and safety tolerance window as function of physical parameters of the nanomagnet including the strength of inertial effects. These results will allow to explain the complex picture arising from numerical simulations of iLLG precessional switching\cite{neeraj2022inertial} where the boundary between successful and unsuccessful switching is extremely and tightly entangled due to a combination of multi-stability and low dissipation that produces the appearance of a quasi-random relaxation process\cite{serpico2009analytical}.

The paper is organized as follows: in section \ref{sec:two}, the iLLG equation is introduced including relationship with relevant material, geometrical parameters of the macrospin particle including the strength of inertial effects the influence of which on precessional switching is the focus of this study. Then, by performing appropriate time rescaling and using suitable systems of coordinates (sections \ref{sec:three},\ref{sec:four}), the analytical approximated solution of magnetization dynamics is developed by using multiple time scale analysis in section \ref{sec:five}. The conditions for realizing successful switching dynamics are addressed in section \ref{sec:six} and, finally, validation of the developed theory with numerical experiments also relevant to experimental investigations is presented in section \ref{sec:seven}. 


\section{Background}\label{sec:two}

We start the analysis by recalling the inertial Landau-Lifshitz-Gilbert (iLLG) equation written for convenience in dimensionless form\cite{DAQUINO2024112874}:
\beq{eq:illg}
\dot{\ue{m}}=\ue{m} \x (\ue{h}_{eff}-\alpha \dot{\ue{m}}-\eta \ddot{\ue{m}}) \mx{,}
\eeq
where $\ue{m}$ is the magnetization unit-vector (normalized by the saturation magnetization $M_s$ of the ferromagnet), time is measured in units of $(\gamma M_s)^{-1}$ ($\gamma$ is the absolute value of the gyromagnetic ratio), $\alpha$ is the Gilbert damping parameter describing dissipative effects (usually $\alpha\sim 10^{-2}$), and the parameter $\eta$ measures the intensity of inertial effects ($\eta=(\gamma M_s \tau)^2$, $\tau$ being the time scale of inertial effects, typically fractions of picoseconds\cite{Neeraj2021inertial} that yields $\eta\sim\alpha \sim 10^{-2}$).  \\
We denote with
\beq{eq:heff}
\ue{h}_{eff}:=\ue{h}_a-\uue{D}\ue{m} \mx{,}
\eeq
the magnetic effective field (also normalized by $M_s$) where $\ue{h}_a \equiv (h_1,h_2,h_3) \in \RR^3$ and $\uue{D}=\diag(D_1,D_2,D_3)$, $D_j \in \RR$ such that 
\beq{eq:easyaxis}
-\infty <D_1 < D_2 < D_3<+\infty \mx{.}
\eeq
Let us recall, see, e.g. \cite{BMS2009}, that the quantity $\Psi(\ue{m}):=|\ue{m}|^2$ is conserved along the solutions of (\ref{eq:illg}), i.e. 
$d \Psi/dt=\ue{m}\cdot \dot{\ue{m}}=0$, 
as $\ue{m} \perp \dot{\ue{m}}$, which is immediately checked from (\ref{eq:illg}), so that any solution starting on $\mathbb{S}^2$ stays on $\mathbb{S}^2$ for all $t \in \RR$. A further standard argument consists in the left cross multiplication of (\ref{eq:illg}) by $\ue{m}$, then using the standard cross product rules and finally the relation $\ue{m}\cdot \ddot{\ue{m}}=-|\dot{\ue{m}}|^2$, in order to write  (\ref{eq:illg}) as follows, 
\beq{eq:illgsixd}
\left\{
\begin{aligned}
    \dot{\ue{m}}&=\ue{v}\\
    \dot{\ue{v}}&= 
-|\ue{v}|^2 \ue{m} -\eta^{-1}\left[\ue{m} \x \ue{v}+ \ue{m} \x (\ue{m} \x \ue{h}_{eff}) + \alpha \ue{v}\right]
\end{aligned}
\right.
\eeq
with $\ue{z}:=(\ue{m},\ue{v}) \in \mathbb{R}^3 \times \mathbb{R}^3=:\ml{Z}$ and for all $\eta \neq 0$. Let us observe that the system of coordinates provided by $\ue{z}$ is conceptually different from the one used in\cite{DAQUINO2024112874}, featuring the variable $\ue{w}:=\ue{m} \times \dot{\ue{m}}$ in place of $\ue{v}$. As the former is more naturally related to the phase space structure of the problem ($\ue{v}$ is simply the velocity), the latter offers a more appropriate geometric interpretation of it and it turns out to be handier for various purposes, such as the geometric integration of the equations proposed in\cite{DAQUINO2024112874}.\\ 
It is immediate to check from (\ref{eq:illgsixd}) that 
\beq{eq:w}
\ml{W}:=\frac{\eta}{2}|\ue{v}|^2+h(\ue{m}), 
\qquad h(\ue{m}):=\frac{1}{2} \left[\sum_{j=1}^3 D_j m_j^2 - D_1\right]-\ue{h}_a \cdot \ue{m}\mx{,}
\eeq
satisfies 
\beq{eq:dwdt}
\frac{d}{dt}\ml{W}=-\alpha |\ue{v}|^2 \mx{,}
\eeq
i.e. $\ml{W}$ is non-increasing along the solutions of (\ref{eq:illgsixd}). \\
Denoting by $\ue{e}_j$ the canonical basis unitary vectors of $\mathbb{R}^3$, let us finally recall that, in the absence of applied field, i.e. $\ue{h}_a=\ue{0}$, the points 
\beq{eq:equilibria}
\ue{z}^{(j),\pm}:=(\pm \ue{e}_j,\ue{0}) \mx{,}
\eeq
are equilibria for the system (\ref{eq:illgsixd}), as it is straightforwardly checked. In particular, a standard analysis via linearisation shows that the points $\ue{z}^{(1),\pm}$ are stable. \\
For this reason, it is meaningful to consider the following class of motions, also known as \emph{switching} 
\begin{defn}\label{def:switching}
The system is said to undergo a switching process if there exists 
a solution $\ue{z}:[0,+\infty) \rw \ml{Z}$ of (\ref{eq:illgsixd}) such that 
\beq{eq:switchingcondition}
\ue{z}(0)=\ue{z}^{(1),+} \qquad ; \qquad 
\lim_{t \rw +\infty}\ue{z}(t)=\ue{z}^{(1),-} \mx{.}
\eeq
\end{defn}
The typical approach to obtaining this class of motions is to apply a suitable field $\ue{h}_a$ for a certain time interval $[0,T_{sw}]$, where $T_{sw}$ is called \emph{time of switching}, then turning it off and wait for the system \emph{relaxation}. Clearly, this requires two main ingredients: $i)$ the attractivity of $\ue{z}^{(1),-}$, i.e., the existence of a non-empty basin of attraction, say $\ml{S}^{-} \ni \ue{z}^{(1),-}$ and $ii)$ a field $\ue{h}_a$ apt to bring the system state in $\ml{S}^{-}$ after a prefixed time $T_{sw}$. As far as $\ue{h}_a$ is concerned, we will restrict ourselves to the case of constant applied fields. As it is reasonable to expect, this operation should ensure a ``certain robustness'' with respect to this instant which, for practical applications, is more properly intended as a window $[T_{sw}-\delta_{sw},T_{sw}+\delta_{sw}]$, commonly termed \emph{safe switching} interval.             


\section{Time rescaled system}\label{sec:three}
The aim of this section is to choose an ``intermediate'' system of coordinates in such a way the subsequent transformation in spherical coordinates of the equations could assume a more convenient structure. To this end, we shall introduce (under suitable assumptions) a new time $\tau$ and a new system of coordinates $\ue{x}$, both defined below. Roughly speaking, the time $\tau$ will allow to ``absorb'' (or, better, ``redirect'') the singularity in $\eta$ of the equations in such a way the new vector field has $O(1)$ magnitude. This is a known approach in the context of singularly perturbed differential equations, see e.g. \cite{Hartmann2010LectureNO}. On the other hand, the coordinates $\ue{x}$ aim to provide a system with respect to which the switching dynamics takes place as a rotation around the $x_3$-axis. Such a system of coordinates will be obtained via a suitably constructed $SO(3)$ transformation.  
\\ For this purpose, let us introduce the following
\begin{hyp}\label{hyp:one}
The term $\ue{h}_a$ shall be supposed to satisfy the following assumptions 
\begin{enumerate}
     \item (Large magnitude) 
\beq{eq:assumptionha}
\ue{h}_a=:\ep^{-1} \hat{\ue{h}}_a \mx{,}
\eeq
where $\ep$ is a ``small'' parameter and $\hat{\ue{h}}_a$ is an object which will be assumed of $O(1)$, i.e. independent on $\ep$. 
\item (``Non-degeneracy''): $h_j >0$ for some $j=2,3$.
\end{enumerate}
\end{hyp}
\begin{rem}\label{rem:ep} Hyp \ref{hyp:one} is justified by the cases of interest in the applications, some mild restrictions will be made in the follow. Experimental studies suggest for $\ep$ values of the order $10^{-1}$. 
\end{rem}
The range for $\ep$ described in Rem. \ref{rem:ep}, compared with the typical values for $\alpha,\eta$ of interest, suggests the following  
\begin{hyp}\label{hyp:two}
We shall assume
\beq{eq:alphaeta}
\alpha=:\ep^2 \hat{\alpha}, \qquad \eta=:\ep^2 \hat{\eta}\mx{,}
\eeq
where, similarly, $\hat{\alpha},\hat{\eta}=O(1)$.
\end{hyp}
\begin{rem} 
It is immediate to realise that the above described choice for $\eta$, suggests that (\ref{eq:illgsixd}) should be treated \emph{de facto} as a singularly perturbed system of ODEs, see e.g. \cite{Hartmann2010LectureNO}. The perturbative analysis of those systems is not elementary but requires some specialised techniques as, for instance, the Multiple Time Scales method we shall employ later on.
\end{rem}
By using Hyp. \ref{hyp:one}, Hyp. \ref{hyp:two} and the expression (\ref{eq:heff}), the iLLG equation (\ref{eq:illg}) reads as
\beq{eq:illg_two}
\ep \dot{\ue{m}}=\hat{\ue{h}}_a \x \ue{m} + \ep \ue{m} \x \uue{D} \ue{m}+\ep^3 \hat{\alpha} \ue{m} \x \dot{\ue{m}} +
\ep^3 \hat{\eta} \ue{m} \x \ddot{\ue{m}} \mx{.}
\eeq
Let us now introduce the following \emph{time rescaling} 
\beq{eq:tau}
t=:\hat{\eta} \ep^2 \tau \mx{,}
\eeq
where $\tau \in \RR$ is the ``new'' time we shall use from now on. By using the latter and denoting with $\ue{f}':=d \ue{f} /d\tau $, for any $\ue{f}=\ue{f}(\tau)$, eq. (\ref{eq:illg_two}) becomes
\beq{eq:illg_three}
\ue{m}'' \x \ue{m} + \ue{m}'=\ep \hat{\eta} \hat{\ue{h}}_a \x \ue{m} + \ep^2 \left[ \hat{\eta} \hat{\alpha} \ue{m} \x \ue{D} \ue{m}+ \hat{\alpha} \ue{m} \x \ue{m}' \right]  \mx{.}
\eeq
It might be useful to notice that the equilibria $\ue{z}^{(j),\pm}$ are mapped to the phase space points $(\ue{m},\ue{m}') \equiv (\pm \ue{e}_j,\ue{0})$ for (\ref{eq:illg_three}) via the transformation (\ref{eq:tau}). \\
Let us now write
\beq{eq:Gamma}
\hat{\ue{h}}_a \x \ue{m}=\ue{\Gamma} \ue{m} \quad ; \quad 
\ue{\Gamma}:= 
\begin{pmatrix}
0&-\hat{h}_3&\hat{h}_2\\
\hat{h}_3&0&-\hat{h}_1\\
-\hat{h}_2&\hat{h}_1&0
\end{pmatrix}
\mx{,}
\eeq
and state the following
\begin{prop}\label{prop:one} Define
\beq{eq:sigmaomega}
\sigma:=\sqrt{\hat{h}_2^2+\hat{h}_3^2}, \qquad \omega:=|\hat{\ue{h}}_a| \mx{,}
\eeq
(note that $\sigma,\omega>0$ by Hyp. \ref{hyp:one}). Then, the matrix
\beq{eq:c} 
\uue{C}:= \frac{1}{\omega}
\begin{pmatrix}
\sigma & 0 & \hat{h}_1\\
-\sigma^{-1} \hat{h}_1 \hat{h}_2 & \sigma^{-1} \omega \hat{h}_3 & \hat{h}_2\\
-\sigma^{-1} \hat{h}_1 \hat{h}_3 & -\sigma^{-1} \omega \hat{h}_2 & \hat{h}_3
\end{pmatrix}
\mx{,}
\eeq
satisfies the following properties
\begin{enumerate}
    \item $\ue{C}\in SO(3)$,
    \item $\uue{C}^T \uue{\Gamma} \uue{C} = \uue{\Lambda}$, 
    where
\beq{eq:Lambda}
 \uue{\Lambda}:=
 \begin{pmatrix}
0&-\omega&0\\
\omega&0&0\\
0&0&0
\end{pmatrix}  \mx{.}
\eeq  
\end{enumerate}
\end{prop}
\proof The proof is postponed to Appendix A. \endproof
Let us now introduce the following coordinates transformation
\beq{eq:x}
\ue{m}:=\uue{C} \ue{x} \mx{.}
\eeq
As a consequence, by Prop. \ref{prop:one}, defining 
\beq{eq:e}
\uue{E}:=\uue{C}^T \uue{D} \uue{C} 
\mx{,}
\eeq
and recalling that $\uue{S}\ue{v} \x \uue{S} \ue{w}=\uue{S}(\ue{u} \x \ue{w})$ for all $\uue{S} \in SO(3)$ and all $\ue{v},\ue{w} \in \RR^3$, it is immediate to check that (\ref{eq:illg_three}) is cast into the following form
\beq{eq:illg_four}
\ue{x}'' \x \ue{x} + \ue{x}'=\ep \hat{\eta} \uue{\Lambda} \ue{x}+ \ep^2 
\left[ \hat{\eta} \ue{x} \x {\ue{E}} \ue{x}+ \hat{\alpha} \ue{x} \x \ue{x}'\right]  \mx{.}
\eeq
\begin{rem} 
The matrix $\uue{E}$ defined in (\ref{eq:e}) is symmetric as $\uue{C} \in SO(3)$ and $\uue{D}$ is diagonal. In particular, its explicit expression (use (\ref{eq:c})), reads as 
\beq{eq:eexplicit}
\uue{E} \equiv 
\begin{pmatrix}
\de \frac{D_1 \sigma^2}{\omega^2}+\frac{\hat{h}_1(D_2 \hat{h}_2^2+D_3 \hat{h}_3^2)}{\sigma^2 \omega^2} & \de \frac{\hat{h}_1 \hat{h}_2 \hat{h}_3 (D_3-D_2)}{\sigma^2 \omega} & 
\de \frac{\hat{h}_1(D_1 \sigma^2-D_2 \hat{h}_2^2-D_3 \hat{h}_3^2)}{\sigma \omega^2} \\[2ex]
\de \frac{\hat{h}_1 \hat{h}_2 \hat{h}_3 (D_3-D_2)}{\sigma^2 \omega} & 
\de \frac{D_2 \hat{h}_3^2+D_3 \hat{h}_2^2}{\sigma^2} & 
\de -\frac{\hat{h}_2 \hat{h}_3(D_3-D_2)}{\sigma \omega}\\[2ex]
\de \frac{\hat{h}_1(D_1 \sigma^2-D_2 \hat{h}_2^2-D_3 \hat{h}_3^2)}{\sigma \omega^2} & 
\de -\frac{\hat{h}_2 \hat{h}_3(D_3-D_2)}{\sigma \omega} & 
\de \frac{D_1 \hat{h}_1^2+D_2 \hat{h}_2^2+D_3 \hat{h}_3^2}{\omega^2}
\end{pmatrix} \mx{.}
\eeq
\end{rem}
\begin{rem}\label{rem:conservationlaw}
As a direct consequence of the fact that $\ue{C}\in SO(3)$, the function $\Psi(\ue{x}):=|\ue{x}|^2$ with $\Psi(\ue{x})\equiv 1$ is conserved along the solutions. On the other hand, this property is immediately checked by recalling $\frac{d}{dt}\Psi=2 \ue{x}\cdot\ue{x}'$ then using (\ref{eq:illg_four}). 
\end{rem}
By defining a new smallness parameter as 
\beq{eq:mu}
\mu:=\sqrt{\hat{\alpha}}\ep
\eeq
and setting 
\beq{eq:defhats}
\hat{\ue{\Lambda}} := 
 \begin{pmatrix}
0&-\hat{\omega}&0\\
\hat{\omega}&0&0\\
0&0&0
\end{pmatrix}, \quad \hat{\omega}:=
\hat{\alpha}^{-\frac{1}{2}} \hat{\eta} \omega
\qquad ; \qquad 
\hat{\uue{E}}:=\hat{\alpha}^{-1} \hat{\eta} \uue{E} \mx{,} 
\eeq 
eq. (\ref{eq:illg_four}) reads as
\beq{eq:illg_final}
\ue{x}'' \x \ue{x} + \ue{x}'=\mu \hat{\ue{\Lambda}} \ue{x}+ \mu^2 
\left[ \ue{x} \x \hat{\uue{E}} \ue{x}+ \ue{x} \x \ue{x}'\right] \mx{,}
\eeq
with initial conditions, see (\ref{eq:switchingcondition}), (\ref{eq:c}) and (\ref{eq:x}), given by
\beq{eq:ic_illg_final}
\ue{x}(0)=\ue{x}^+:= \uue{C}^{T} \ue{e}^{(1)} = \omega^{-1} (\sigma,0,- \hat{h}_1 )\qquad ; \qquad \ue{x}'(0)=\ue{0} \mx{.}
\eeq
Let us notice that the target point we aim to reach has coordinates
\beq{eq:xminus}
\ue{x}^-:= \uue{C}^{T} \left(-\ue{e}^{(1)} \right) = \omega^{-1} (-\sigma,0, \hat{h}_1)\qquad \mx{.}
\eeq


\section{Dynamics in spherical coordinates}\label{sec:four}
As stressed in Rem. \ref{rem:conservationlaw}, the dynamics in the variables $\ue{x}$ inherits the property to take place on the unit sphere. For this purpose, let us consider the standard transformation $\mathrm{S}:(r,\theta,\ph) \rw \ue{x}$ defined by 
\beq{eq:spherical}
\ue{x}=(r \sin \theta \cos \ph, r \sin \theta \sin \ph,r \cos \theta ) 
\eeq
and denote with $\hat{e}_{i,j}$ the entries of $\hat{\uue{E}}$. For all $\delta>0$, let us define the \emph{poles-less sphere} as
\beq{eq:polesless}
\mathbb{S}_{\delta}^2:=[0,2 \pi] \times [\delta, \pi-\delta] \ni (\ph,\theta) \mx{,}
\eeq
then state the following
\begin{prop}\label{prop:equivalence} Set $\ml{T}>0$. Then  (\ref{eq:illg_final}) can be written in the coordinates system given by (\ref{eq:spherical}) and reads as
\beq{eq:illg_spherical_lhs}
\left\{
\begin{aligned}
\sin \theta \ph''+2 \cos \theta \ph' \theta' + \theta' &=\mu^2 \ml{F}_1 \\
\theta''- \cos \theta \sin \theta (\ph')^2- \sin \theta \ph' +\mu \hat{\omega} \sin \theta &=\mu^2 \ml{F}_2 
\end{aligned}
\right.
\mx{,}
\eeq
where 
\beq{eq:illg_spherical_rhs}
\begin{aligned}
\ml{F}_1 &= \left[(\hat{e}_{1,1}-\hat{e}_{2,2}) \cos \ph \sin \ph - 2  \hat{e}_{1,2} \cos^2 \ph + \hat{e}_{1,2}- \ph' \right] \sin \theta \\
& + \left[ \hat{e}_{1,3} \sin \ph - \hat{e}_{2,3} \cos \ph \right] \cos \theta \\
\ml{F}_2 &= 
\left[(\hat{e}_{2,2}-\hat{e}_{1,1}) \cos^2 \ph - 2 \hat{e}_{1,2} \cos \ph \sin \ph + \hat{e}_{3,3}- \hat{e}_{2,2} \right] \sin \theta \cos \theta \\
& + \left[ \hat{e}_{2,3} \sin \ph + \hat{e}_{1,3} \cos \ph  \right] (2 \sin^2 \theta -1) - \theta' 
\end{aligned} 
\eeq
provided that $(\ph(\tau),\theta(\tau)) \in \mathbb{S}_{\delta_{\ml{T}}}^2$ for all $\tau\in[0,\ml{T}]$ and some  $\delta_{\ml{T}}>0$.
\end{prop}
\proof
If $(\ph(t),\theta(t)) \in \mathbb{S}_{\delta_T}^2$, then $\sin \theta \neq 0$, so that the Jacobian matrix of the transformation $\pl \ue{x}/\pl (r,\ph,\theta)$ is invertible. Hence, it is sufficient to substitute (\ref{eq:spherical}) in   (\ref{eq:illg_final}) and finally left-multiply both sides by $\left(\pl \ue{x}/\pl (r,\ph,\theta)\right)^{-1}$ to obtain (\ref{eq:illg_spherical_lhs}) and (\ref{eq:illg_spherical_rhs}).  
\endproof
Let us now consider the coordinates implicitly defined by the following  
\beq{eq:chixi}
\ph=: \mu \hat{\omega} \tau + \mu \xi \quad ; \quad \theta=:\pi/2+\mu \chi \mx{.}
\eeq
By substituting the latter into (\ref{eq:illg_spherical_lhs}) and (\ref{eq:illg_spherical_rhs}), then dividing by $\mu$, we obtain
\beq{eq:illg_chixi}
\left\{
\begin{aligned}
\cos(\mu \chi) \xi'' - 2 \mu \sin(\mu \chi)(\hat{\omega}+\xi')\chi'+\chi' &=\mu \tilde{\ml{F}}_1 \\
\chi''+2^{-1} \mu \sin (2 \mu \chi) (\hat{\omega}+ \chi')^2 - \cos(\mu \chi) \xi' &=\mu \tilde{\ml{F}}_2 
\end{aligned}
\right.
\mx{,}
\eeq
where 
\beq{eq:illg_chixi_rhs}
\begin{aligned}
\tilde{\ml{F}}_1 &= [(\hat{e}_{1,1}-\hat{e}_{2,2})\cos(\mu \chi)\cos(\mu \hat{\omega} \tau +\mu \xi)-\hat{e}_{1,3}\cos(\mu \chi)]\sin(\mu \hat{\omega} \tau +\mu \xi)\\
& + \hat{e}_{2,3} \sin(\mu \chi) \cos(\mu  \hat{\omega} \tau+\mu \xi)
- 2 \hat{e}_{1,2} \cos(\mu \chi(\tau)) \cos^2(\mu  \hat{\omega} \tau+\mu \xi)\\
& + (\hat{e}_{1,2}-\mu \xi') \cos(\mu \chi(\tau))\\
\tilde{\ml{F}}_2 &=[2 \mu \hat{e}_{1,2} \cos(\mu \chi) \sin(\mu \chi)  \cos(\mu  \hat{\omega} \tau+\mu \xi) +2 \hat{e}_{2,3} \cos^2 (\mu \chi)-\hat{e}_{2,3}]\sin(\mu \hat{\omega} \tau +\mu \xi)\\
& +  \hat{e}_{1,3}(2 \cos^2(\mu \chi)-1)\cos(\mu  \hat{\omega} \tau+\mu \xi) + \mu (\hat{e}_{1,1}-\hat{e}_{2,2}) \cos(\mu \chi) \sin(\mu \chi)  \cos^2(\mu  \hat{\omega} \tau+\mu \xi)\\
&+\mu (\hat{e}_{2,2}-\hat{e}_{3,3}) \cos(\mu \chi) \sin(\mu \chi) - \mu \chi' \mx{.}
\end{aligned}
\eeq
denote the functions $\ml{F}_j$ with $j=1,2$, computed in the new coordinates $(\chi,\xi)$. \\
Let us now recall that a system of the form 
\beq{eq:illg_integrable}
\left\{
\begin{aligned}
\xi''  + \chi' &=f_1(\tau) \\
\chi'' - \xi' &=f_2(\tau) 
\end{aligned}
\right.
\mx{,}
\eeq
can be explicitly solved via elementary methods. More precisely, by setting 
\beq{eq:z}
\ue{q}\equiv (q_1,q_2,q_3,q_4):=(\xi,\xi',\chi,\chi')\mx{,} 
\eeq
(\ref{eq:illg_integrable}) can be written as 
\beq{eq:illg_approx_fo}
\ue{q}'= \uue{A} \ue{q}+\mu \ue{\mathfrak{f}}(\tau) \mx{,}
\eeq
where 
\[
\uue{A}=
\begin{pmatrix}
0 & 1& 0 & 0\\
0 & 0 &0 &-1\\
0 & 0 &0 &1\\
0 & 1 & 0 & 0
\end{pmatrix}, \quad
\ue{\mathfrak{f}}(\tau)=
\begin{pmatrix}
0\\
f_1(\tau)\\
0\\
f_2(\tau)
\end{pmatrix}
\mx{,}
\]
so that the solution reads as 
\beq{eq:variationconstants}
\ue{q}(\tau)=e^{\uue{A} \tau}\ue{q}(0)+\int_0^{\tau} e^{\uue{A} (s-\tau)} \ue{\mathfrak{f}}(s) ds \mx{,} \qquad 
e^{\uue{A}\tau}=
\begin{pmatrix}
0 & \sin \tau & 0 & \cos \tau -1\\
0 & \cos \tau &0 &-\sin \tau\\
0 & 1- \cos \tau &0 &\sin \tau\\
0 & \sin \tau & 0 & \cos \tau
\end{pmatrix} \mx{.}
\eeq
\begin{rem}
A direct Taylor expansion in $\mu$ of the circular functions appearing in (\ref{eq:illg_chixi}) and (\ref{eq:illg_chixi_rhs}) immediately yields, up to $O(\mu)$, a system of the form (\ref{eq:illg_integrable}), where $\ue{f}(\tau) \equiv \mu (-\hat{e}_{1,2} ,\hat{e}_{1,3})$. \\
However, such an expansion is meaningful not only if the functions $\chi,\xi$ stay bounded throughout the evolution, namely $|\chi(\tau)|,|\xi(\tau)| \leq \ml{K}$, with $\ml{K}=O(1)$, but also if $\tau$ does. Although the former is reasonably assumed (more on this later), the latter turns out to be a too restrictive condition in our case. More precisely, as we will be interested in $O(\mu^{-1})$ timescales, terms associated with  $\sin(\mu  \hat{\omega} \tau+\cdot)$ are no longer negligible.
\end{rem}
This difficulty is naturally dealt with the method of Multiple Time Scales, see, e.g. \cite{nay}, as shown in the next section.


\section{Approximated solution: the method of Multiple Time Scales}\label{sec:five}

\begin{lem}\label{lem:approximation} The solution to (\ref{eq:illg_spherical_lhs}) is approximated, up to $O(\mu^2)$, by
\beq{eq:thphapp}
\begin{aligned}
\ph^{[\leq 2]}(\tau)&=\ph(0)+\mu \hat{\omega} \tau + 
\mu \left[\hat{\omega}^{-1} \left(\hat{e}_{1,3} \sin (\mu \hat{\omega} \tau)-
\hat{e}_{2,3} \cos (\mu\hat{\omega} \tau) + \hat{e}_{2,3}\right)-\hat{\omega} \sin(\tau) \right]\\
& + \mu^2 \left[\hat{e}_{1,2} (1-\cos \tau) - \hat{e}_{1,3} \sin \tau \right],\\     
\theta^{[\leq 2]}(\tau)&=\theta(0) + \mu \left[-\hat{\omega} (1-\cos \tau)+(2 \hat{\omega})^{-1}  \left((\hat{e}_{2,2}-\hat{e}_{1,1})\sin^2(\mu \hat{\omega} \tau)-\hat{e}_{1,2}\sin(2 \mu \hat{\omega} \tau)\right) \right]\\
&+ \mu^2 \left[(\hat{e}_{2,3}-\hat{e}_{1,3})(1-\cos \tau)-\hat{e}_{1,2} \sin \tau - \hat{e}_{2,3} \sin (\mu \hat{\omega} \tau) \cos \tau \right]
\end{aligned}
\eeq
i.e. $\ph-\ph^{[\leq 2]}(\tau),\theta(\tau)-\theta^{[\leq 2]}(\tau)=O(\mu^3)$.\\
Furthermore, for all $\theta(0)\in (0,\pi)$, define $\tilde{d}:=\min\{\theta(0),\pi-\theta(0)\}$ and
\beq{eq:muzero}
\mu_0:=\left[2 \hat{\omega} + |\hat{e}_{1,2}|(1+1/(2 \hat{\omega})) + 
2 |\hat{e}_{1,3}| + 3 |\hat{e}_{2,3}|+(|\hat{e}_{1,1}|+|\hat{e}_{2,2}|)/(2 \hat{\omega})\right]^{-1} \tilde{d}
\mx{.}
\eeq
Then, for any $\mu \in(0,\mu_0]$, one has $(\ph^{[\leq 2]}(\tau),\theta^{[\leq 2]}(\tau)) \in \mathbb{S}_{\tilde{d}/2}^2$ for all $\tau \in \mathbb{R}$, so that (\ref{eq:thphapp}) provide an approximation for the solution to (\ref{eq:illg_final}) via the mapping (\ref{eq:spherical}) (with $r \equiv 1$).
\end{lem}
\begin{rem} The statement implies that, at the price of a possible further restriction of $\mu$ with respect to the threshold $\mu_0$, the full (exact) solution satisfies $(\ph(\tau),\theta(\tau)) \in \mathbb{S}_{\tilde{r}}^2$ with $\tilde{r}<\tilde{d}$, i.e. it does never reach the poles of $\mathbb{S}^2$. 
\end{rem}
\proof
Following \cite{nay} we set 
\beq{eq:times}
T_0:=\tau, \qquad T_1:=\mu \tau \mx{,}
\eeq
which yield \cite{nay}, 
\beq{eq:rules}
d/d\tau=\pl_{T_0}+\mu \pl_{T_1} + O(\mu^2), \qquad
d^2/d\tau^2=\pl_{T_0}^2+2 \mu \pl_{T_0} \pl_{T_1} + O(\mu^2) \mx{,}
\eeq
and we seek solutions to (\ref{eq:illg_chixi}) of the form 
\beq{eq:solutions}
\begin{aligned}
\chi & =\chi(T_0,T_1)=\chi^{(0)}(T_0,T_1)+\mu \chi^{(1)}(T_0,T_1)+O(\mu^2),\\ 
\xi & =\xi(T_0,T_1)=\xi^{(0)}(T_0,T_1)+\mu \xi^{(0)}(T_0,T_1)+O(\mu^2) \mx{.} 
\end{aligned}
\eeq
It is now sufficient to substitute the expansions (\ref{eq:solutions}) into (\ref{eq:illg_chixi}) using the rules (\ref{eq:rules}) for the derivatives involved, to get the following systems
\beq{eq:illg_mult_zero}
\left\{
\begin{aligned}
\pl_{T_0}^2 \xi^{(0)}  + \pl_{T_0} \chi^{(0)} &=0 \\
\pl_{T_0}^2 \chi^{(0)} - \pl_{T_0} \xi^{(0)} &=0 
\end{aligned}
\right.
\mx{,}
\eeq
and 
\beq{eq:illg_mult_one}
\left\{
\begin{aligned}
\pl_{T_0}^2 \xi^{(1)}  + \pl_{T_0} \chi^{(1)} &= \pl_{T_1} \chi^{(0)} + 2 \pl_{T_0 T_1}^2 \xi^{(0)} + (\hat{e}_{1,1}-\hat{e}_{2,2}) \sin (T_1 \hat{\omega}) \cos (T_1 \hat{\omega})\\
&+(1-2 \cos^2 (T_1 \hat{\omega})) \hat{e}_{1,2} \\
\pl_{T_0}^2 \chi^{(1)} - \pl_{T_0} \xi^{(1)} &= - \pl_{T_1} \xi^{(0)} + 2 \pl_{T_0 T_1}^2 \chi{(0)} + \hat{e}_{2,3} \sin (T_1 \hat{\omega})+ \hat{e}_{1,3} \cos (T_1 \hat{\omega})  
\end{aligned}
\right.
\mx{,}
\eeq
obtained equating the powers of $\mu^0$ and $\mu^1$, respectively.\\
It is important to notice that the systems (\ref{eq:illg_mult_zero}) and, in particular (\ref{eq:illg_mult_one}), are both of the form (\ref{eq:illg_integrable}). Hence, they can be resolved explicitly. It is sufficient to change the definition (\ref{eq:z}) into $\ue{z}_j\equiv (z_1^{(j)},z_2^{(j)},z_3^{(j)},z_4^{(j)}):=(\xi^{(j)},\pl_{T_0} \xi^{(j)},\chi^{(j)},\pl_{T_0} \chi^{(j)})$ for $j=0,1$, then use (\ref{eq:variationconstants}). However, as it is customary in the Multiple Time Scales method, the object $\ue{z}_0(0)$ is no longer a constant but will be a function of $T_1$. More precisely, by redefining it as $\ue{w}(T_1)$, the solution to (\ref{eq:illg_mult_zero}) reads as   
\beq{eq:solutionlevelzero}
\begin{aligned}
\xi_0(T_0,T_1)&=\sin(T_0) w_2(T_1)-\left(1-\cos(T_0) \right)w_4(T_1)+w_1(T_1)\\
\chi_0(T_0,T_1)&=\sin(T_0) w_4(T_1)+\left(1-\cos(T_0) \right)w_2(T_1)+w_3(T_1)
\end{aligned}
\mx{,}
\eeq
where $\ue{w}(T_1)$ is a function which will be chosen later in order to annihilate potential secular terms. \\
Although with these extra ``degrees of freedom'', the zero-th order solution can be thought as completely determined. Hence, by substituting (\ref{eq:solutionlevelzero}), it is now immediate to observe that the system (\ref{eq:illg_mult_one}) is, once more, of the form (\ref{eq:illg_integrable}), as anticipated. The only difference consists in the presence of a non-homogeneous term, which is straightforwardly determined by computing the extra-terms appearing in (\ref{eq:illg_mult_one}) explicitly, i.e.  
\beq{eq:nonhomogeneousterm}
\begin{aligned}
f_1(T_0,T_1)&=(\hat{e}_{1,1}-\hat{e}_{2,2}) \cos( \hat{\omega} T_1) \sin( \hat{\omega} T_1) + \hat{e}_{1,2}\left(1-2 \cos^2 ( \hat{\omega} T_1) \right)\\
&+ 2 \cos(T_0) w_2' -\sin(T_0) w_4'+
\left(1-\cos (T_0) \right) w_2' w_3'\\ 
f_2(T_0,T_1)&=\hat{e}_{1,3} \cos( \hat{\omega} T_1)+ \hat{e}_{2,3} \sin( \hat{\omega} T_1)\\
&+ \left(1+ \cos (T_0) \right) w_4' + \sin (T_0) w_2'-w_1'
\end{aligned}
\mx{,}
\eeq
where $w_j':=d/d T_1 w_j$. It is now possible to solve (\ref{eq:illg_mult_one}) via (\ref{eq:variationconstants}). However, this is exactly the point in which one realises that if $\ue{w}_j(T_1)$ are not supposed to be functions, but simple constants, secular terms (i.e. linear in $T_0$) would inevitably arise. The obtained expression are slightly cumbersome so they will be omitted, but these offer a guideline to craft the suitable vector $\ue{w}_j(T_1)$ during this process of removal of the mentioned unwanted terms. The suitability of the vector constructed in this way
\beq{eq:choicew}
\begin{aligned}
w_1^*&:=\hat{\omega}^{-1}
\left[ 
\hat{e}_{1,3} \sin (\hat{\omega} T_1) - 
\hat{e}_{2,3} \sin (\hat{\omega} T_1)
\right] + C_1 \\
w_2^*&:=C_2\\
w_3^*&:=(2 \hat{\omega})^{-1}
\left[ 
(\hat{e}_{1,1}-\hat{e}_{2,2}) \cos^2 (\hat{\omega} T_1)+ \hat{e}_{1,2} \sin (2 \hat{\omega} T_1) \right]+C_3 \\
w_4^*&:=C_4
\end{aligned}
\eeq
with $C_j \in \mathbb{R}$ to be determined, is straightforwardly checked to be the ``right one'' as it has the effect to bring the contribution of the whole non-homogeneous term $f_{1,2}$ to \emph{zero}, i.e. $f_{1,2}|_{\ue{w}=\ue{w}^*}=0$, so that the solution of (\ref{eq:illg_mult_one}) is the trivial one
\beq{eq:solutionpostw}
\xi_1 = \chi_1 \equiv 0, \qquad \forall T_0 \in \mathbb{R} \mx{.} 
\eeq
As a consequence, by (\ref{eq:solutions}), the required solution up to the first order in $\mu$, is immediately obtained from the expressions (\ref{eq:solutionlevelzero}) for $\xi_0,\chi_0$, by substituting the newfound (\ref{eq:choicew}) in place of $\ue{w}$. This yields  
\beq{eq:finalxichi}
\begin{aligned}
    \xi&=\hat{\omega}^{-1}
\left[ \hat{e}_{1,3} \sin (\hat{\omega} T_1) - \hat{e}_{2,3} \sin (\hat{\omega} T_1) \right] + C_2 \sin(T_0) + C_4 (\cos(T_0)-1) + C_1\\
    \chi&=(2 \hat{\omega})^{-1}
\left[ 
(\hat{e}_{1,1}-\hat{e}_{2,2}) \cos^2 (\hat{\omega} T_1)+ \hat{e}_{1,2} \sin (2 \hat{\omega} T_1) \right]+ 
C_4 \sin(T_0) + C_2 (1-\cos(T_0)) + C_3 \\
&+ \mu \hat{e}_{2,3}(1-\cos(T_0))(1+ \sin (\hat{\omega} T_1))
\end{aligned}
\mx{.}
\eeq
It is now sufficient to replace the two timescales as functions of $\tau$ (\ref{eq:times}) then impose the initial condition. This leads to the determination of the constants $C_j$  
\beq{eq:constantsc}
\begin{array}{rclrcl}
C_1&=&\chi(0)+\hat{\omega}^{-1}\hat{e}_{2,3}, & C_2&=&\xi'(0)-\mu \hat{e}_{1,3}, \\
C_3&=&\chi(0)+(2 \hat{\omega})^{-1} (\hat{e}_{2,2}-\hat{e}_{1,1}), & C_4&=&\chi'(0)-\mu \hat{e}_{1,2} \mx{.}
\end{array} 
\eeq
The last step consists in going back to the original angular coordinates. For this purpose, let us recall formulae (\ref{eq:chixi}), which allow to write the initial conditions $\chi(0)$ and $\xi(0)$ appearing in (\ref{eq:constantsc}) in terms of those for $\ph(0)$ and $\theta(0)$. More precisely 
\[
\xi(0)=\ph(0)/\mu, \quad \xi'(0)=-\hat{\omega}, \quad 
\chi(0)=(\theta(0)-\pi/2)/\mu, \quad \chi'(0)=0 \mx{.} 
\]
By using these relations in (\ref{eq:constantsc}) and employing (\ref{eq:chixi}) once more, we obtain the desired expressions (\ref{eq:thphapp}).\\
As for the threshold $\mu_0$, this is straightforwardly obtained from (\ref{eq:thphapp}) by using that $\mu<1$ by hypothesis, so that $\theta^{[\leq 2]}(\tau)-\theta(0)$ is bounded by $\mu$ times the upper bound of the coefficients of both $\mu$ and $\mu^2$.
\endproof

\section{Switching dynamics}\label{sec:six}
The aim of this section is to establish criteria for the existence of a proper switching dynamics. As already stressed, a successful switching necessarily relies on a certain attractivity property of the target equilibrium. The latter is established by the following 
\begin{lem}[Attractivity]\label{lem:attractivity} Let us consider  (\ref{eq:illg}) in the form (\ref{eq:illgsixd}),  set  
\beq{eq:dji}
D_{j,i}:=D_j-D_j, \qquad 1 \leq i \leq j \leq 3 \mx{,}
\eeq
note that $D_{j,i}>0$ by (\ref{eq:easyaxis}), then define
\beq{eq:basin}
\mathcal{S}:=\{(\ue{m},\ue{v}) \in \ml{Z} \, : \, \ml{W}(\ue{m},\ue{v}) \leq D_{2,1}/3\} \mx{.}
\eeq
Clearly $\ue{z}^{(1),\pm} \in \mathcal{S}$ by (\ref{eq:dwdt}) and (\ref{eq:easyaxis}). Let us now denote with $\mathcal{S}^{-}$ and $\mathcal{S}^{+}$ the connected components of $\mathcal{S}$ containing $\ue{z}^-$ and $\ue{z}^+$, respectively. \\
Then, any initial condition for (\ref{eq:illgsixd}) starting in $\mathcal{S}^-$ ($\mathcal{S}^+$) is attracted by the (asymptotically stable) equilibrium $\ue{z}^-$ ($\ue{z}^+$). In other terms, $\mathcal{S}^{\pm}$ are contained in the basin of attraction of $\ue{z}^{\pm}$, respectively. 
\end{lem}
\begin{rem}
As it is common in Lyapounov-type arguments and Stability Theory in general, the one given by $\mathcal{S}^{\pm}$, might only be a (small) part of the whole basins of attraction. This is a very common issue, see, e.g. \cite{khal}, \cite{fliap}. It is very well-known, on the other hand, that (more) refined estimates of such a region could potentially require \emph{ad-hoc} and highly specific analyses. 
\end{rem}
\proof We shall discuss the case of $\ml{S}^{-}$, being $\ml{S}^{+}$ analogous, \emph{mutatis mutandis}. Let us recall (\ref{eq:w}) and the first order setting (\ref{eq:illgsixd}) for (\ref{eq:illg}). It is easy to show, via a standard constrained Optimisation argument, that the function $h(\ue{m})$ possesses on $\mathbb{S}^2$ the following stationary points 
\beq{eq:hej}
h(\pm \ue{e}_j) = D_{j,1}/2, \quad j=1,2,3 \mx{.}
\eeq
In particular, $h(\pm \ue{e}_1)=0$, hence, by (\ref{eq:easyaxis}), $\pm \ue{e}_1$ are the only minima. This implies that  the set of all the points $\ue{m} \in \mathbb{S}^2$ such that $h(\ue{m})< D_{2,1}/2$ contains minima only and, in particular, $-\ue{e}_1$ is the only minimum in the connected component of this set containing $-\ue{e}_1$ itself. As a consequence, the set $\mathcal{S}^-$ contains \emph{a fortiori} the only equilibrium point $\ue{z}^{(1),-}$ for (\ref{eq:illgsixd}). This set is clearly compact by definition, but it is also positively invariant, being the sub-levels set of the non-increasing function $\ml{W}$. \\
Furthermore, as $\alpha \neq 0$ by assumption,  (\ref{eq:dwdt}) implies that $d \ml{W}/dt = 0$ if, and only if, $\ue{v}=0$. However, from (\ref{eq:illgsixd}), the only point of $\mathcal{S}^-$ satisfying this condition is the equilibrium $\ue{z}^{(1),-}$ itself. The proof is now reduced to a direct consequence of the Barbashin - Krasovskii - La Salle's Theorem (see, e.g. \cite{lasalle}).   
\endproof
Let us first introduce the following 
\begin{hyp}\label{hyp:third}
The applied field $\ue{h}_a$ satisfies
\beq{eq:thirdhyp}
h_1 h_3=0 \mx{,}
\eeq
i.e., $\hat{h}_1 \hat{h}_3=0$, by (\ref{eq:assumptionha}).
\end{hyp}
This assumption has the purpose to simplify the implementation of the definition of $T_{sw}$, avoiding, at the same time, technicalities of marginal mathematical interest. 
It is worth observing that it does not represent a significant restriction with respect to fields $\ue{h}_a$ commonly used in the applications. In fact, typical prototypes have the form $\ue{h}_a=(0,h_2,0)$ or $\ue{h}_a=(0,0,h_3)$, although variants of the first case with $h_1 \neq 0$ may be considered as well. Hypotheses \ref{hyp:one} and \ref{hyp:third} are satisfied in the mentioned cases. By using the material presented in the previous sections we can finally state the main result as follows  
\begin{satz}[Switching]\label{thm:main} Suppose $\hat{\ue{h}}_a$ of the form described in Hyp \ref{hyp:one}, 2. and that its components, as well as $\hat{\alpha}$ and $\hat{\eta}$, are given positive $O(1)$ constants. \\
Let us now compute $\sigma, \omega$ by (\ref{eq:sigmaomega}), then $\hat{\omega}$ and $\hat{e}_{i,j}$ via (\ref{eq:defhats}) and  (\ref{eq:eexplicit}), and finally $\mu_0$ according to (\ref{eq:muzero}).\\
Then, for any sufficiently small $\mu \in (0,\mu_0]$ and any initial condition $\dot{\ue{m}}(0)=\ue{0}$ and $\ue{m}(0) \in \mathbb{S}^2$ for which $\mathrm{S}^{-1} \left(\ue{C}^{-1}\ue{m}(0) \right)$ satisfies $\theta(0)\neq 0,\pi$, the solution $\ue{m}(t)$ of Eq. (\ref{eq:illg}) in which $\alpha$ and $\eta$ have been determined with (\ref{eq:mu}) and (\ref{eq:alphaeta}), is approximated by $\ue{m}^{[\leq 2]}(\tau):=\ue{C}  \mathrm{S} (1,\theta^{[\leq 2]}(\tau),\ph^{[\leq 2]}(\tau))$ up to $O(\mu^2)$. Similarly, the solution of (\ref{eq:illgsixd}) is approximated up to $O(\mu)$ (see proof for the explicit expressions). \\
Let us now include Hyp. \ref{hyp:third} and define 
\beq{eq:newmu}
\tilde{\mu}_0:=2 \pi [4 \hat{\omega}^2 + |\hat{e}_{1,1}|+|\hat{e}_{2,2}|]^{-1} \hat{\omega} \mx{,}
\qquad \Xi:=(8 \hat{\alpha} \hat{\omega}^2)^{-1}D_{2,1}\hat{\eta}
\mx{.}
\eeq
Then, for all sufficiently small $\mu \in (0,\tilde{\mu}_0]$, and any applied field of the form
\beq{eq:switchoff}
\ue{h}_a(t):=
\left\{ 
\begin{aligned}
\ue{h}_a &\qquad  t \in [0,t^*) \\
0 & \qquad t \geq t^* 
\end{aligned}
\right.
\eeq
where $t^* \in [T_{sw} - \delta_{sw},  T_{sw} + \delta_{sw}]$ (safe switching interval), 
\beq{eq:tdelta}
T_{sw}:=\mu \frac{  \hat{\eta}}{\hat{\omega} \hat{\alpha}} \pi, \quad 
\delta_{sw}= \mu\frac{ \hat{\eta}}{\hat{\omega} \hat{\alpha}} \arcsin
\left(
\min\left\{
1,
\sqrt{\mu^2 \hat{\omega}^2 + K_f} - \mu \hat{\omega}\right\} 
\right) 
\eeq
and 
\beq{eq:Kf}
K_f:=D_{3,1}/(4 D_{2,1}) \mx{,}
\eeq
then the system undergoes a successful switching according to Def. \ref{def:switching} if 
\begin{itemize}
\item[i] $\Xi \geq 1$, or
\item[ii] $\hat{\omega}=1/(2 \mu n)$, for some $n \in \mathbb{N}$ and $t^* \in [T_{sw} - \delta_{sw}^*,  T_{sw} + \delta_{sw}^*]$ where 
\beq{eq:deltaswstar}
\delta_{sw}^*:=\min\{\delta_{sw}, \mu^2 (\hat{\eta}/\hat{\alpha})\arccos(1-\Xi) \} \mx{.}
\eeq
\end{itemize}
\end{satz}
\begin{rem} The new threshold $\tilde{\mu}_0$ is nothing but the $\mu_0$ defined in (\ref{eq:muzero}) after the simplification introduced by Hyp. \ref{hyp:third}. It is interesting to observe that the condition on $\mu \leq \tilde{\mu}_0$ is necessary to use this statement but it is not sufficient. The experiments of the next section will give practical example of when  $\mu$ might or might not be considered ``sufficiently small'' .\\
A more ``operative'' version of this theorem is presented in Appendix B in algorithm form.
\end{rem}
This statement has two main consequences. The first part ensures that if $\mu$ is sufficiently small \emph{and} not exceeding $\mu_0$, then the system dynamics can be approximated via explicit formulae. The second part is focused on the switching motion, hence in the particular case of the initial condition $\ue{z}(0) \equiv \ue{z}^{(1),+}$. In this hypothesis, under the additional assumption described in Hyp. \ref{hyp:third}, the results guarantees the existence of a switching solution in two possible scenarios. The first one occurs  when $\Xi \geq 1$. The latter identifies a full manifold of parameters: recalling (\ref{eq:defhats}), this is equivalent to the condition $8 \hat{\eta} |\hat{\ue{h}}_a|\leq D_{2,1} $. However, if the latter is not satisfied (which might be the case of some applications), $ii)$ provides a ``backup case'', although this is clearly valid only on a discrete set of values for $\hat{\omega}$ (indexed by $n$) and at the price of a potential restriction of the safe switching interval. The case study of the next section will provide an example of this ``dual'' validity this main statement.
\proof[Proof of Thm. \ref{thm:main}]
As $\theta(0) \neq 0,\pi$ by assumption, the quantity $\tilde{d}$ of Lemma \ref{lem:approximation} is strictly positive and so is, as a consequence, the threshold $\mu_0$ as defined in (\ref{eq:muzero}). Hence, the first part of the statement is a direct consequence of Lemma \ref{lem:approximation}.\\
With the aim to focus on the switching dynamics, let us now introduce Hyp. \ref{hyp:third}. The latter leads to a remarkable simplification of the approximated equations structure. More precisely, as it implies 
\beq{eq:simplye}
\hat{e}_{1,2}=\hat{e}_{1,3}=\hat{e}_{2,3}=0 \mx{,}
\eeq
(recall that $\uue{E}$ is symmetric), the $SO(3)$ transformation induced by $\uue{C}$ simply reduces to the permutation 
\beq{eq:cnew}
 \uue{C}:=
 \begin{pmatrix}
1 & 0 &0\\
0 & 0& 1\\
0&-1&0
\end{pmatrix}  \mx{,}
\eeq  
implying that 
\beq{eq:newe}
\uue{E}=\diag(D_1,D_3,D_2) \mx{.}
\eeq
Furthermore, by Def. \ref{def:switching}, the initial condition is of the form described by the first of (\ref{eq:switchingcondition}), which implies, in particular, $\ph(0)=0$ and $\theta(0)=\pi/2$. In the notation of Lemma \ref{lem:approximation}, this implies $\tilde{d}=\pi/2$. As anticipated, let us now denote with $\tilde{\mu}_0$, the threshold $\mu_0$ of formula (\ref{eq:muzero}). The latter, because of (\ref{eq:simplye}), has exactly the form (\ref{eq:newmu}). \\      
As a consequence, if $\mu \in (0,\tilde{\mu}_0]$ we can use Lemma \ref{lem:approximation} once again so that, starting from (\ref{eq:thphapp}), using (\ref{eq:simplye}), and the transformations (\ref{eq:spherical}) and (\ref{eq:x}) afterwards, the approximated equations can be readily cast into the original set of variables
\beq{eq:apporiginal}
\begin{aligned}
m_1^{[\leq 2]}(\tau)&=
\cos \left[ \mu (2 \hat{\omega})^{-1}(\hat{e}_{2,2}-\hat{e}_{1,1})\sin^2 (\mu \hat{\omega} \tau)- \mu \hat{\omega} (1-\cos \tau) \right]
\cos \left[ \mu \hat{\omega}(\sin \tau - \tau)\right] 
\\
m_2^{[\leq 2]}(\tau)&=-
\sin \left[ \mu (2 \hat{\omega})^{-1}(\hat{e}_{2,2}-\hat{e}_{1,1})\sin^2 (\mu \hat{\omega} \tau)- \mu \hat{\omega} (1-\cos \tau) \right]
\\
m_3^{[\leq 2]}(\tau)&=
\cos \left[ \mu (2 \hat{\omega})^{-1}(\hat{e}_{2,2}-\hat{e}_{1,1})\sin^2 (\mu \hat{\omega} \tau)- \mu \hat{\omega} (1-\cos \tau) \right]
\sin \left[ \mu \hat{\omega}(\sin \tau - \tau)\right] 
\end{aligned}
\eeq
Clearly, it is immediate to check that $\ue{m}^{[\leq 2]}(0)=\ue{e}_1$ holds as initial condition. Furthermore, it is immediate to check that 
by setting 
\beq{eq:tausw}
\ml{T}_{sw}:=(\mu \hat{\omega})^{-1} \pi \mx{,}
\eeq
one has
\beq{eq:expressionattsw}
\begin{aligned}
\ue{m}^{[\leq 2]}(\ml{T}_{sw})&=\left(
-\cos[\mu \hat{\omega} (\cos \ml{T}_{sw}-1)] 
\cos[\mu \hat{\omega} \sin \ml{T}_{sw}], \right. \\
&-\sin[\mu \hat{\omega} (\cos \ml{T}_{sw}-1)],\\
& \left.-\cos[\mu \hat{\omega} (\cos \ml{T}_{sw}-1)] 
\sin[\mu \hat{\omega} \sin \ml{T}_{sw}] \right) \mx{,}
\end{aligned}
\eeq
which is clearly $O(\mu)$ close to $-\ue{e}_1$. On the other hand, the value (\ref{eq:tausw}) is immediately derived from the coordinates (\ref{eq:chixi}) as the time $\tau$ such that $|\ph(0)-\ph(\tau)|=\pi$.\\
As we are interested in the behaviour of these solutions for sufficiently small $\mu$, a first order expansion in $\mu$ would suffice, provided that the latter is non-trivial. This implies imply that, within the $O(\mu^{-1})$ timescale represented by $\ml{T}_{sw}$, the solution will be accurate up to $O(\mu)$. Hence, by defining $\ue{m}^{[\leq 1]}$ the first order Taylor expansion in $\mu$ of $\ue{m}^{[\leq 2]}$, it is easy to verify that  
\beq{eq:apporiginalone}
\begin{aligned}
m_1^{[\leq 1]}(\tau)&=
\cos (\mu \hat{\omega} \tau) + \mu \sin (\mu \hat{\omega} \tau) \sin \tau  
\\
m_2^{[\leq 1]}(\tau)&=
- \mu \left[ (2 \hat{\omega})^{-1}(\hat{e}_{2,2}-\hat{e}_{1,1})\sin^2 (\mu \hat{\omega} \tau)- \hat{\omega} (1-\cos \tau) \right]
\\
m_3^{[\leq 1]}(\tau)&=
-\sin (\mu \hat{\omega} \tau) + \mu \cos (\mu \hat{\omega} \tau) \sin \tau  
\end{aligned}
\eeq
The full description of an approximate solution for (\ref{eq:illgsixd}), clearly involves the derivatives of (\ref{eq:apporiginalone}) with respect to $\tau$ , i.e. $\ue{z}^{[\leq 1]}:=(\ue{m}^{[\leq 1]},d \ue{m}^{[\leq 1]}/d\tau)$. However, this has to be done quite carefully as a straightforward expansion of trigonometric functions depending on the ``slow'' time $\mu \hat{\omega} \tau$ would cease to be accurate over $O(\mu^{-1})$ timescales. A workaround consists in reintroducing once more the timescales (\ref{eq:times}), i.e. $\tilde{\ue{m}}^{[\leq 1]}:=\ue{m}^{[\leq 1]}|_{(\mu \tau=T_1,\tau=T_0)}$, computing the time derivative of $\tilde{\ue{m}}^{[\leq 1]}$ via the first of the rules (\ref{eq:rules}), then Taylor-expanding it at the first order in $\mu$ and finally restoring the time $\tau$ via (\ref{eq:times}). This leads to the following expressions, which, by definition, is up to $O(\mu)$  
\beq{eq:mdot}
d \ue{m}^{[\leq 1]}/d\tau \equiv \mu \hat{\omega} ( \sin (\mu \hat{\omega} \tau) (\cos \tau-1), 
 \sin \tau,
 \cos (\mu \hat{\omega} \tau) (\cos \tau-1)) \mx{.}
\eeq
The expression (\ref{eq:apporiginalone}) will be used as approximation for the (exact) $\ue{m}(\tau)$ (which correspond to the $\ue{m}(t)$ of Sec. \ref{sec:two} via the time transformation (\ref{eq:tau})), in order to evaluate the value attained by them at the instant $\mathcal{T}_{sw}$ and establish whether or not the trajectory reaches the basin of attraction $\mathcal{S}^{-}$ provided by Lemma \ref{lem:attractivity}. Clearly, this discussion could be carried out in two equivalent ways: change the time in $\mathcal{W}$ first, then use (\ref{eq:apporiginalone}) and (\ref{eq:mdot}) or change the time in (\ref{eq:apporiginalone}) and (\ref{eq:mdot}) then substitute in the original $\mathcal{W}$ afterwards. This would imply a proof in the original time $t$ or in the transformed time $\tau$, respectively. The first approach will be used here. \\
For this purpose, let us recall that by (\ref{eq:tau}) and (\ref{eq:mu}), the formula expressing the transformation between the two time variables reads as
\beq{eq:newtau}
t=\mu^2 (\hat{\eta}/\hat{\alpha}) \tau \mx{.}
\eeq
Let us now denote with $\hat{\mathcal{W}}$ the function $\mathcal{W}$ written in the time-rescaled variables $\ue{m}(\tau)$. The only difference consists in the term of $\mathcal{W}$ containing the derivative with respect to $t$, which has to be transformed by using (\ref{eq:newtau}). In this way, we obtain
\beq{eq:what}
\hat{\mathcal{W}}(\ue{m}(\tau))=\mu^{-2} \frac{\hat{\alpha}}{2 \hat{\eta}} |\ue{m}'(\tau)|^2+h(\ue{m}(\tau)) \mx{.}
\eeq
If we require that 
\beq{eq:conditionbasin}
\hat{\mathcal{W}} \left(\ue{m}^{[\leq 1]} (\ml{T}_{sw}) \right) \leq D_{2,1}/4 \mx{,}
\eeq
then, for sufficiently small $\mu$, one necessarily has $\hat{\mathcal{W}}(\ue{m}(\ml{T}_{sw}) \leq  D_{2,1}/3$ i.e. the $O(\mu^2)$ remainder satisfies 
\[
\left| \hat{\mathcal{W}}(\ue{m}(\ml{T}_{sw})- \hat{\mathcal{W}} \left(\ue{m}^{[\leq 1]} (\ml{T}_{sw}) \right) \right|\leq D_{2,1}/3-D_{2,1}/4=D_{2,1}/12 \mx{.} 
\]
The proof relies on the fact that if (\ref{eq:conditionbasin}) holds with strict inequality, then such a property persists in a whole neighbourhood of the switching time $\tau=\ml{T}_{sw}$ by continuity, giving rise to the safe switching interval described in the statement.\\
In order to achieve this task, let us start showing that 
\beq{eq:secondaryineq}
h \left(\ue{m}^{[\leq 1]} (\tau) \right) \leq  D_{2,1}/8 \mx{,}
\eeq
for all $\tau$ in a suitable interval (safe switching). For this purpose, let us introduce in (\ref{eq:apporiginalone}) and (\ref{eq:mdot}) the time translation
\beq{eq:s}
\tau=:\ml{T}_{sw}+s \mx{,}
\eeq
define $\tilde{\ue{m}}^{[\leq 1]}(s):=\ue{m}^{[\leq 1]}(\ml{T}_{sw}+s)$, then evaluate $h(\tilde{\ue{m}}^{[\leq 1]}(s))$ in a neighbourhood of $s=0$. A Taylor expansion in $\mu$ shows that
\beq{eq:thirdineq}
h(\tilde{\ue{m}}^{[\leq 1]}(s)) =  D_{3,1} \left[ 2^{-1}\sin^2 (\mu \hat{\omega} s) - \mu \hat{\omega} \sin(\mu \hat{\omega} s) \cos(\mu \hat{\omega} s) \sin ((\mu \hat{\omega} s + \pi)/(\mu \hat{\omega})) \right]  + O(\mu^2) \mx{.}
\eeq
If $O(\mu^2)$ are disregarded, the r.h.s. of the latter vanishes for $s=0$, so that (\ref{eq:secondaryineq}) clearly holds for $\tau=\ml{T}_{sw}$. A proof for (\ref{eq:tdelta}) easily follows by observing that a bound of the r.h.s. of (\ref{eq:thirdineq}) is given by $2^{-1}D_{3,1}(\sin^2 (\mu \hat{\omega} s)+2 \mu \hat{\omega} |\sin (\mu \hat{\omega} s)|)$. Hence, the instants for which the latter is less or equal than $D_{2,1}/8$ are given by the positive solutions of the elementary inequality $y^2+2\mu \hat{\omega}y\leq K_f$, recall (\ref{eq:Kf}), via the variable transformation $y:=|\sin (\mu \hat{\omega} s)|$. The second of (\ref{eq:tdelta}) is easily recovered by coming back to the original time $t$. \\
The proof is complete if a similar bound is shown for the first term appearing in $\hat{\mathcal{W}} \left(\ue{m}^{[\leq 1]} (\ml{T}_{sw}) \right)$ as well. For this purpose, let us notice that the latter, by (\ref{eq:mdot}), reads as 
\beq{eq:inequalitykin}
(2 \mu^2 \hat{\eta})^{-1} \hat{\alpha} |d \ue{m}^{[\leq 1]}(\tau)/d \tau|^2 =\hat{\eta}^{-1} \hat{\alpha} \hat{\omega}^2 (1-\cos \tau) \leq D_{2,1}/8 \mx{.}
\eeq
If $\Xi \geq 1$ as in case $i)$, the last inequality of (\ref{eq:inequalitykin}) clearly holds (uniformly in $\tau$). Otherwise $\Xi<1$, thus $\delta_{sw}^*$ in (\ref{eq:deltaswstar}) is correctly defined. Moreover, if $\hat{\omega}$ satisfies the condition described in case $ii)$, then $\ml{T}_{sw} \equiv 2 n \pi$ so that $(1-\cos \tau)$ vanishes at $\tau=\ml{T}_{sw}$. This means that (\ref{eq:inequalitykin}) holds in a whole neighbourhood of $\tau=\ml{T}_{sw}$, provided that $\tau \in [\ml{T}_{sw}-\arccos(1-\Xi), \ml{T}_{sw}+\arccos(1-\Xi)]$. The latter immediately gives the interval appearing in case $ii)$ via (\ref{eq:newtau}). This completes the proof. \\
It is also interesting to observe from (\ref{eq:expressionattsw}) that, in this case $ii)$, one gets $\ue{m}^{[\leq 1]}(\ml{T}_{sw})=-\ue{e}_1$ exactly.
\endproof


\section{Validation and numerical experiments}\label{sec:seven}
The aim of this section is to provide a numerical validation of Thm. \ref{thm:main}. This task will be addressed in two different cases. In the first case (sec. \ref{subsec:one}), more ``theoretically-oriented'', the attention is focused on the capability of the main result in producing an approximation of the solutions, showing that its quality improves as $\mu$ is reduced. For this reason, a ``paradigmatic'' choice of parameters is considered without worrying too much about their experimental significance. In the second case (sec \ref{subsec:two}), an experimentally relevant set of parameters is considered. Incidentally, the latter constitutes a ``stress test'' for Thm. \ref{thm:main} as the values required for $\mu$ are more than ten times higher w.r.t. the first case study. Nevertheless, the most important approximation properties (for instance of the first component of the magnetisation), can still be appreciated, as well as the possibility to ensure a successful switching. 

\subsection{Case study $1$: ``Theoretical''}\label{subsec:one}
\begin{figure*}[t!]\begin{center}
  		\begin{minipage}[c][1\width]{0.49\textwidth}
			\hspace{0pt}
			{\begin{overpic}[width=1\textwidth]{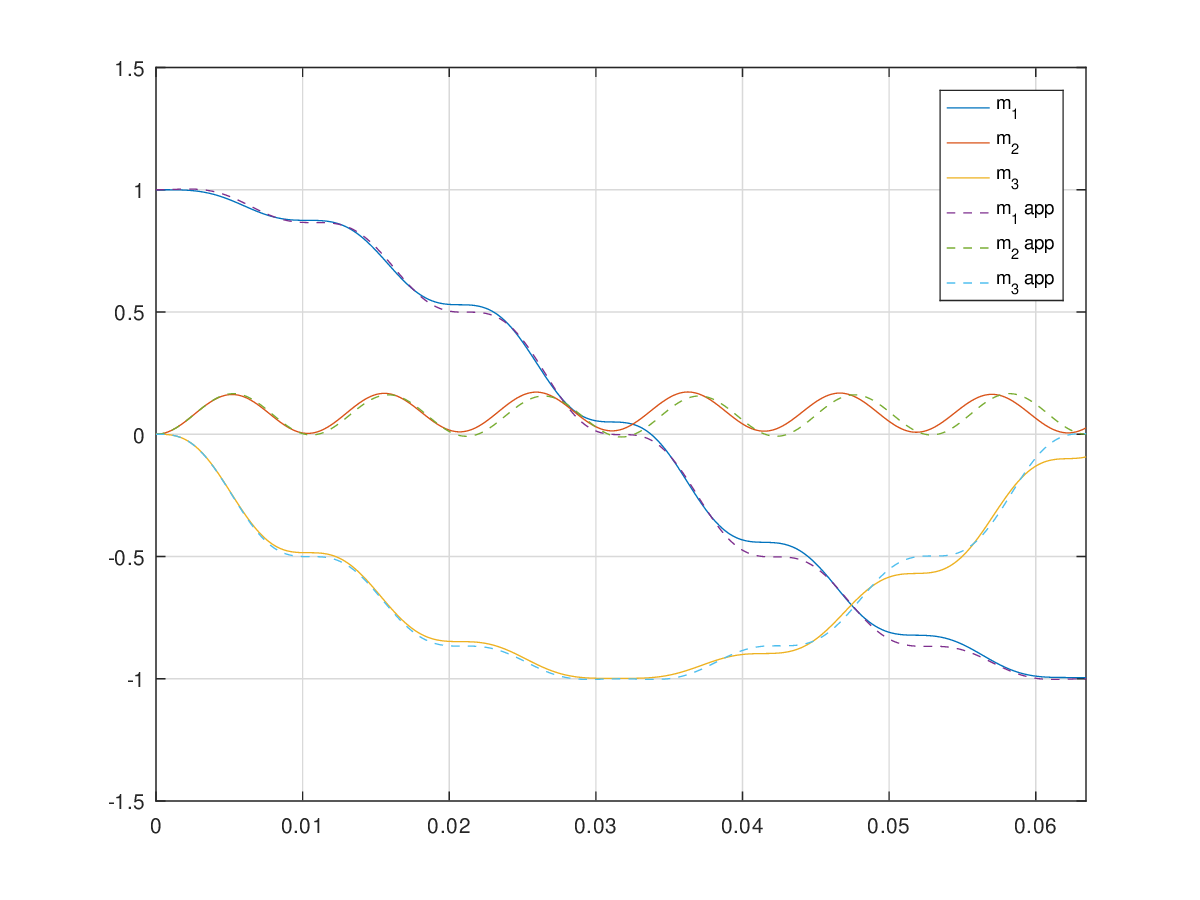}
					\put(0,78){\footnotesize (a)}
			\end{overpic}}
		\end{minipage}	
		\begin{minipage}[c][1\width]{0.49\textwidth}
			\hspace{5pt} 
			{\begin{overpic}[width=\textwidth]{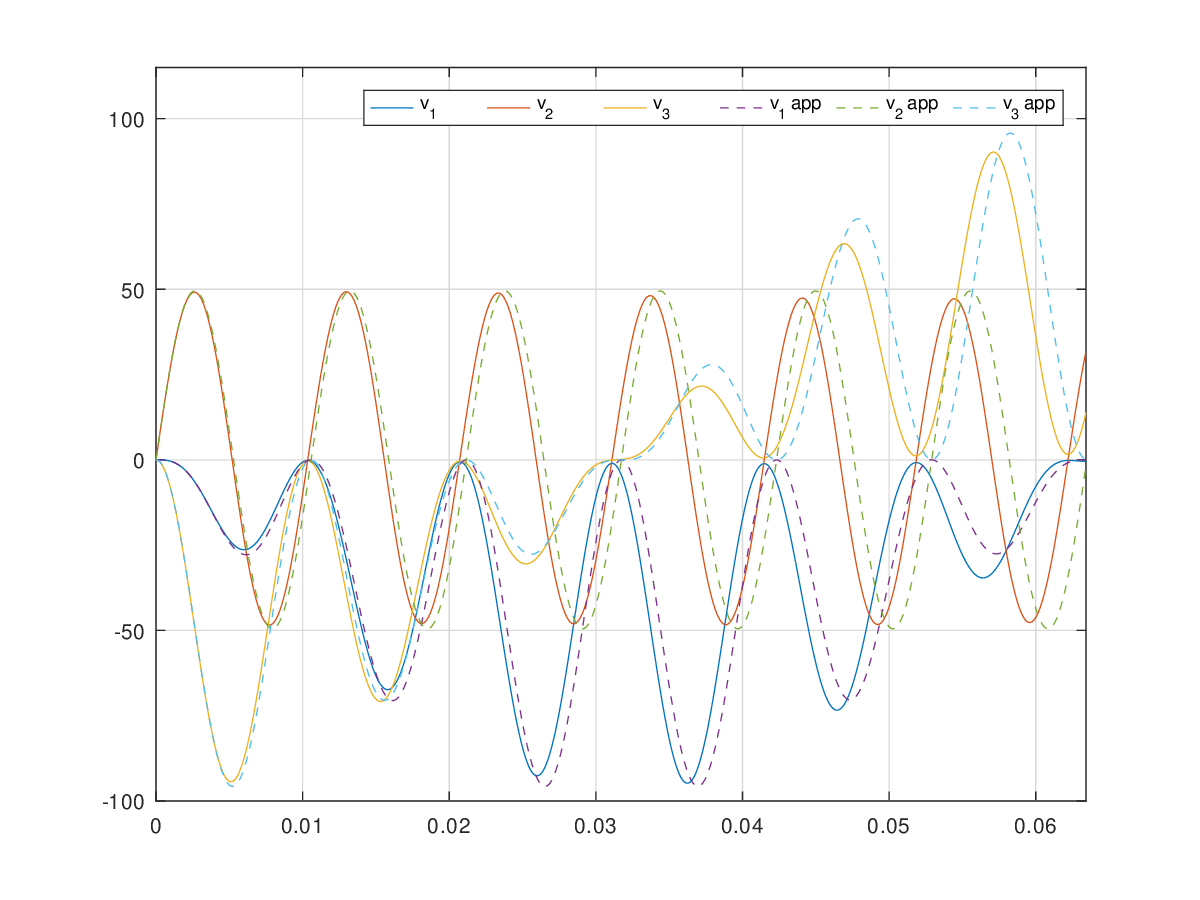}
					\put(0,78){\footnotesize (b)}
			\end{overpic}}
		\end{minipage}\\
		\vspace{-30pt}
  	\begin{minipage}[c][1\width]{0.49\textwidth}
			\hspace{5pt} 
			{\begin{overpic}[width=\textwidth]{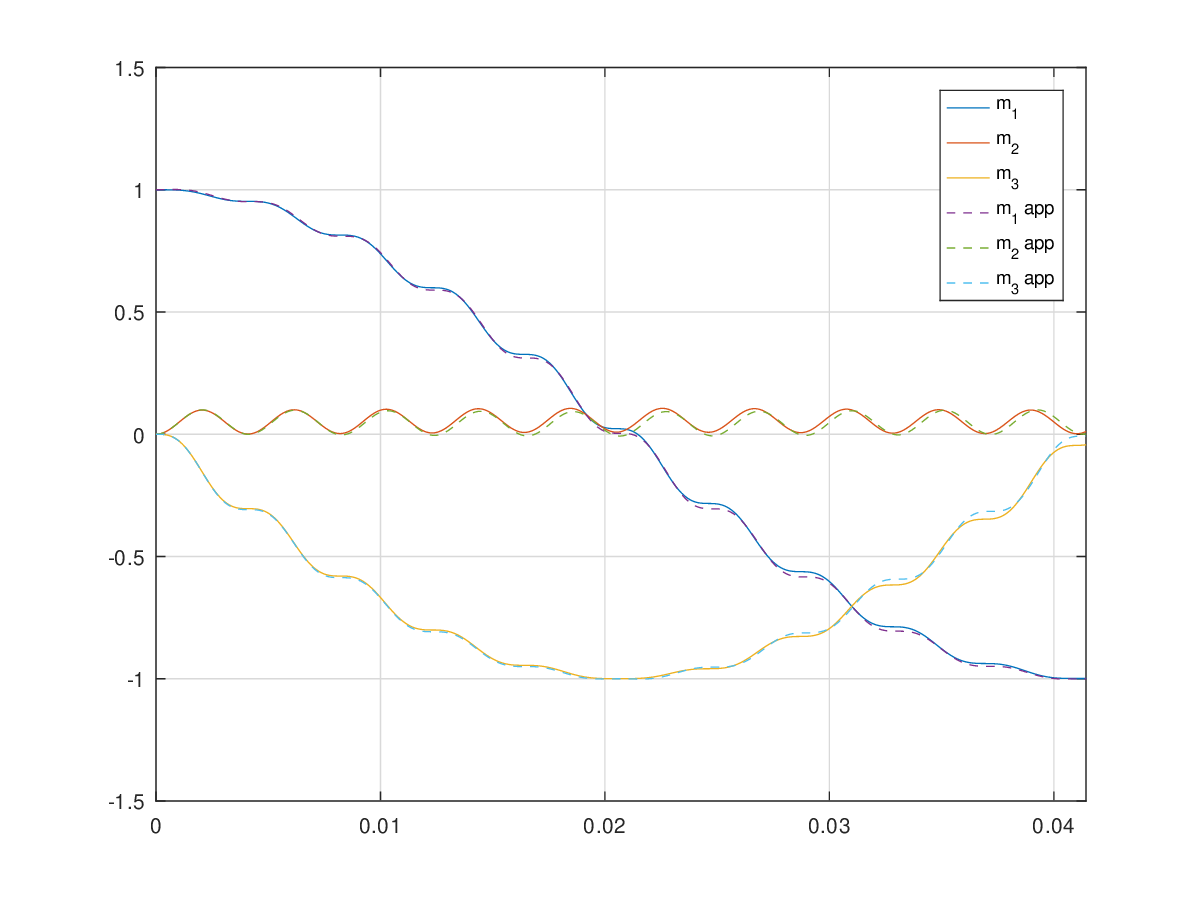}
					\put(0,78){\footnotesize (c)}
			\end{overpic}}
		\end{minipage}	
    \begin{minipage}[c][1\width]{0.49\textwidth}
			\hspace{5pt}
			{\begin{overpic}[width=\textwidth]{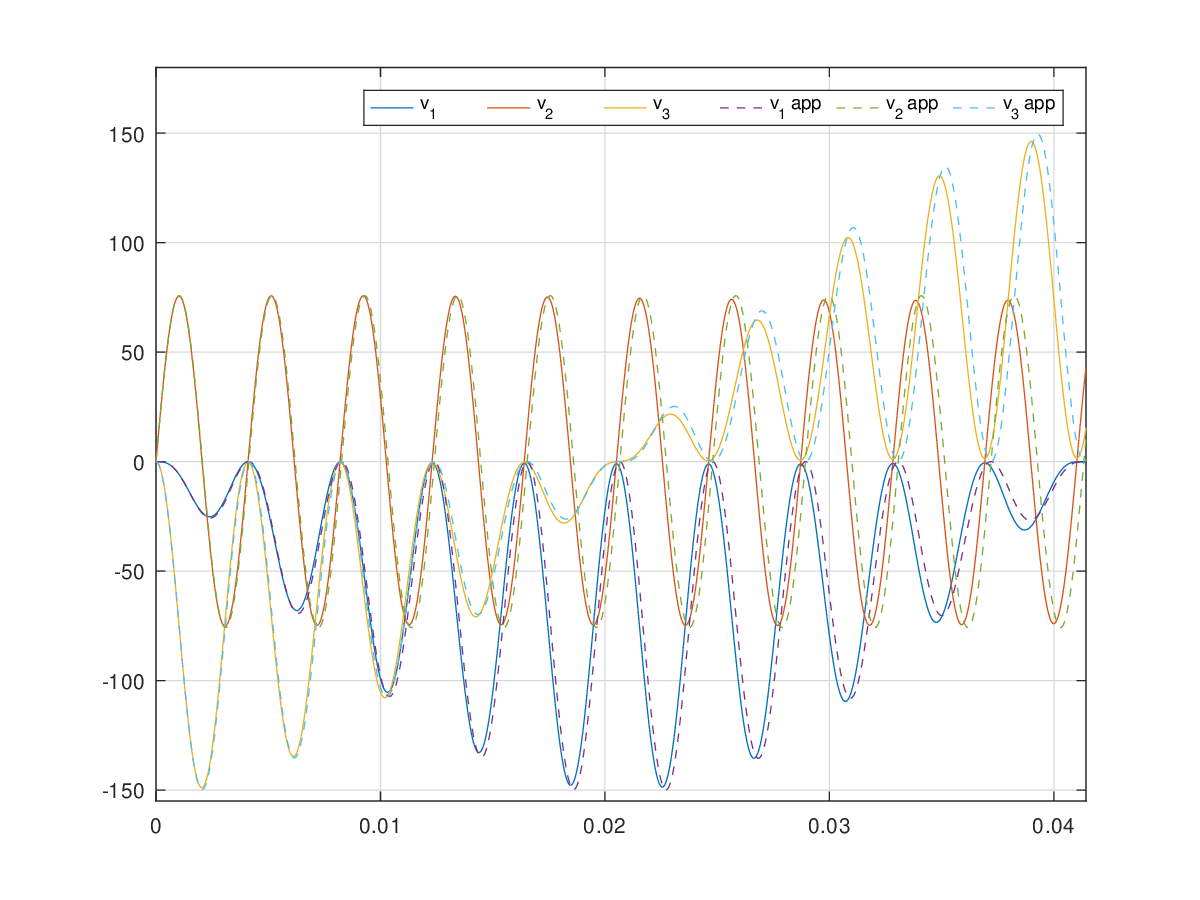}
					\put(0,78){\footnotesize (d)}
 			\end{overpic}}
		\end{minipage}\\
	\end{center}
	\vspace{-30pt}
	\caption{Comparison between the numerically computed solutions (continuous line) of (\ref{eq:illgsixd}) and their approximation provided by (\ref{eq:apporiginalone}) and (\ref{eq:mdot}) (dashed line) for $\mu_A$ in panels (a), (b) and $\mu_B$ in panels (c), (d), respectively. Let us recall that (\ref{eq:apporiginalone}) and (\ref{eq:mdot}) are meant to be plotted via the time transformation (\ref{eq:tau}). Panels (a) and (c) show the components of $\ue{m}(t)$ whilst the corresponding derivatives $\ue{v}(t)$ are shown in panels (b) and (d).}
	\label{fig:validation1}
	\vspace{0pt}
\end{figure*}

\begin{figure*}[t!]\begin{center}
		\begin{minipage}[c][1\width]{0.48\textwidth}
			\hspace{0pt} 
			{\begin{overpic}[width=\textwidth]{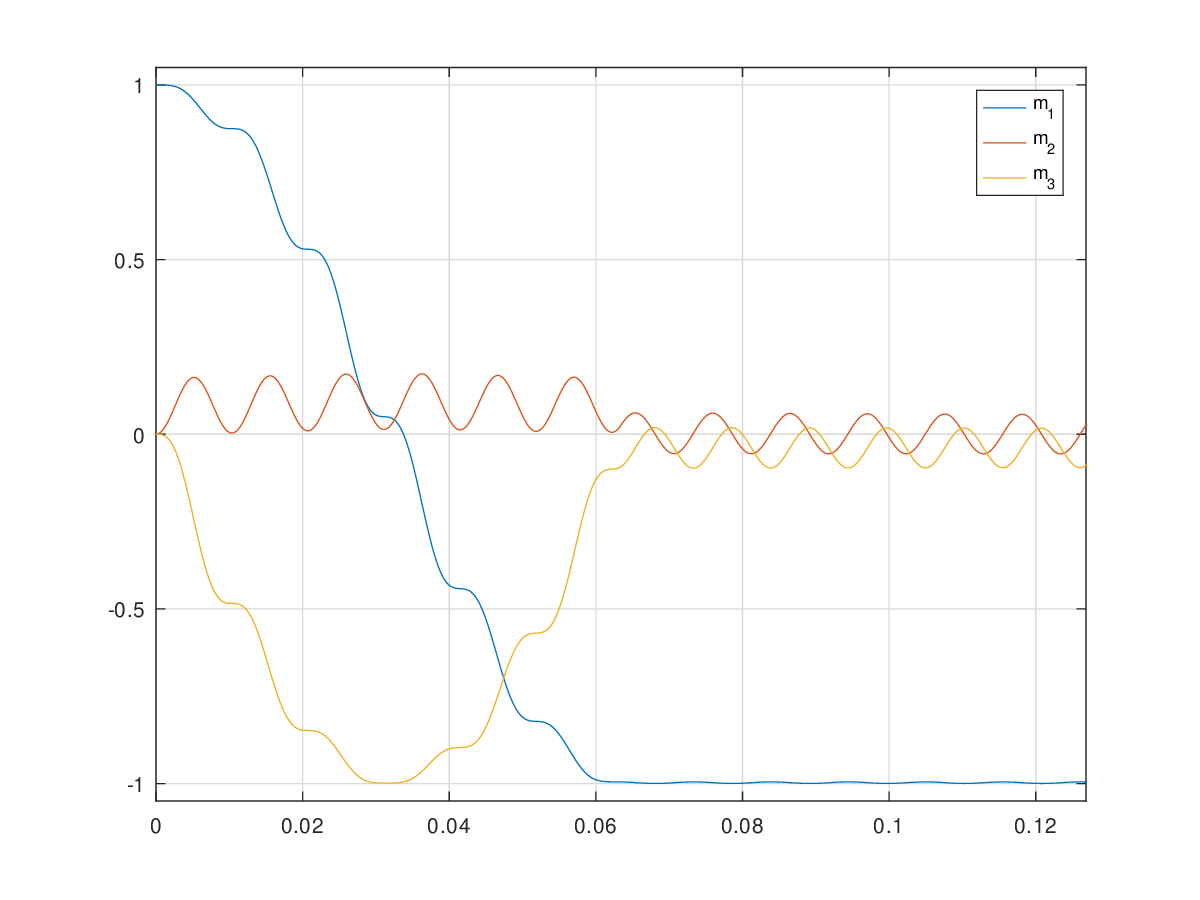}
					\put(0,78){\footnotesize (a)}
					\put(68,2){\footnotesize $t$}
                \end{overpic}}
		\end{minipage}
    		\begin{minipage}[c][1\width]{0.48\textwidth}
			\hspace{-5pt}
			{\begin{overpic}[width=1\textwidth]{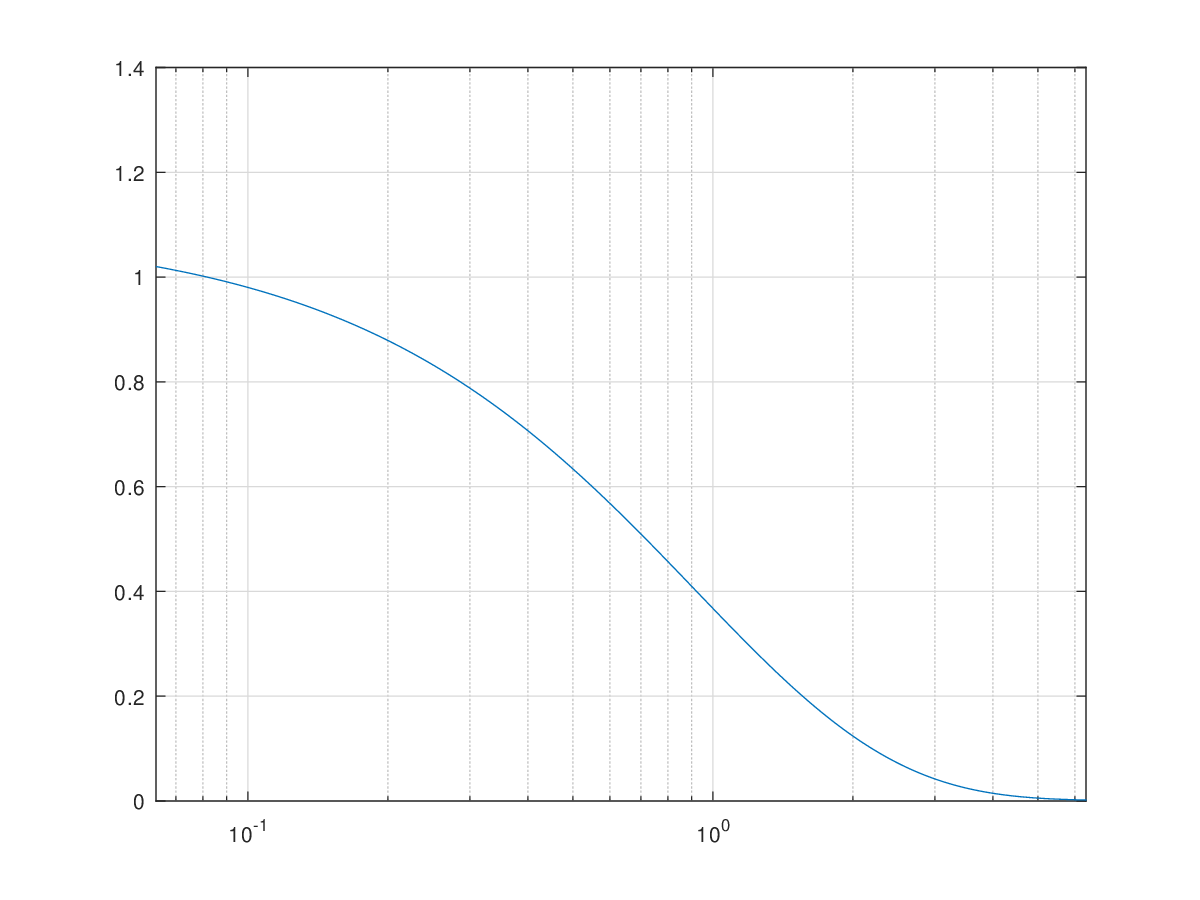}
					\put(0,78){\footnotesize (b)}
					\put(68,2){\footnotesize $t$}
			\end{overpic}}
		\end{minipage}	\\
		\vspace{-30pt}
  \end{center}
	\caption{In panel (a), the system evolution for $t \in [0, 2 T_{sw}]$, where $T_{sw}=0.0635$ is the switch-off instant of the applied field according to (\ref{eq:switchoff}) with $t^*=T_{sw}$. Panel (b) shows the decay of the function $\mathcal{W}$ during the system relaxation for $t \in [T_{sw},100 T_{sw}]$ ($\log$ scale). }
	\label{fig:validation2}
	\vspace{0pt}
\end{figure*}

To this end, let us consider a model described by the following parameters
\begin{itemize}
    \item $\uue{D}=(-0.1087,0,1)$, effective demagnetizing factors of a thin-film with in-plane anisotropy, see sec.\ref{subsec:two},
    \item $\hat{\ue{h}}_a=(0,b,0)$ where $b$ is an $O(1)$ constant to be determined, 
    \item $\hat{\alpha}=2.3$, $\hat{\eta}=4.21$.
\end{itemize}
In order to proceed with the validation we are following the statement of Thm. \ref{thm:main} closely or, equivalently, the algorithms reported in Appendix B. However, as Hyp. \ref{hyp:third} is also satisfied, the Algorithm \ref{alg:two} can be used, although the order of some steps needs to be rearranged as $\hat{\ue{h}}_a$ is not fixed but depends upon the parameter $b$. We find 
\beq{eq:sigmahatomega}
\sigma=\omega=b, \qquad  \hat{\omega}=2.7760 b \mx{.}
\eeq
By (\ref{eq:newe}), we have $\uue{E}=(-0.1087,1,0)$ hence  $\hat{\uue{E}}=(-0.1990,1.8304,0)$ by (\ref{eq:defhats}). \\
From (\ref{eq:newmu}) we find, in particular 
\[
\Xi=0.003227/ b^2 \mx{.}
\]
so that we are in case $i)$ for all $b \leq 0.05681$. With the aim to choose $\mu \sim 10^{-2}$ this value for $b$ is not in agreement with the requirement for it to be $O(1)$ (w.r.t. $\mu$). Hence, we shall investigate the options given by $ii)$. According to this case, we can only consider a finite set of values for $b$, and more precisely all those of the form 
\beq{eq:b}
b_n=(2 \mu \hat{\eta} n)^{-1}\sqrt{\hat{\alpha}}=0.1801/(\mu n) \mx{.} 
\eeq
For instance, if we consider a ``sufficiently small'' $\mu_A=0.03033$ (corresponding to $\ep=0.0200$), we have from (\ref{eq:b})
\beq{eq:bn}
b_n^{(A)}=5.9439/n\mx{,}
\eeq
hence, the only admissible values according to $ii)$, i.e. such that $\Xi_n:=\Xi \left( b_n^{(A)} \right)=9.1339 \cdot 10^{-5} n^2 <1$ are obtained for $n=1,\ldots,104$. In fact $\Xi_{104}=0.9879$ whilst $\Xi_{105}=1.0070$. In the light of this admissible range and of (\ref{eq:bn}), we shall set $b$ as $ b_6^{(A)} =0.9906$.\\
On the other hand, $\tilde{\mu}_0(b)=17.4421 b/(30.8245 b+2.0294)$, hence $\tilde{\mu}_0\left( b_6^{(A)} \right)=0.5353$ which is fully compatible with the requirement $\mu_A \leq \tilde{\mu}_0$. \\
Let us now choose a smaller value of $\mu$, e.g. $\mu_B=0.0190$, in this case we obtain $b_n^{(B)}=9.4789/n$. It is easy to check that a wider range for $n$ is allowed is this case, and more precisely $n\leq 166$. We can choose, for instance, $n=10$ in such a way to obtain a similar $b=0.9479$. The inequality $\mu_B \leq \tilde{\mu}_0=0.5562$ is amply satisfied as well. The problem parameters which correspond to the two choices for $\mu$ are reported in tab. \ref{tab:one}
\begin{table}[H]
\begin{center}
\begin{tabular}{|c|l|l|l|l|l|l|}
\hline
$\mu$  & $\alpha$ & $\eta$ & $\omega$ & $(\ue{h}_a)_2$ \\
\hline
	$\mu_A$  & 9.20 $\cdot$ $10^{-4}$ & 1.68 $\cdot$ $10^{-3}$ & 0.9906  & 49.53\\
\hline
	$\mu_B$  & 3.59 $\cdot$ $10^{-4}$ & 6.58 $\cdot$ $10^{-4}$ & 0.9479  & 78.83\\
\hline
\end{tabular}
\caption{Values of the parameters which correspond to the above mentioned choices for $\mu$. Here $(\ue{h}_a)_2$ denotes the second (i.e. the non-zero) component of $\ue{h}_a$.}
\label{tab:one}
\end{center}
\end{table}
Fig. \ref{fig:validation1} shows the validation of the approximation formulae (\ref{eq:apporiginalone}) and (\ref{eq:mdot}) derived in the proof of Thm. \ref{thm:main} for the two described cases of $\mu_A$ and $\mu_B$, via a comparison with the numerical solution of system (\ref{eq:illgsixd}) with the initial conditions described in Def. \ref{def:switching}. Let us recall that these expression constitute a particular case of the ``general'' approximation Lemma \ref{lem:approximation}.\\
An example of the possibility to obtain a complete switching phenomenon as described in Thm. \ref{thm:main}, is shown in Fig. \ref{fig:validation2}. In this case, obtained by choosing once more $\mu=0.0303$ and $n=6$ in (\ref{eq:bn}), the applied field $\ue{h}_a$ is switched off at $t=T_{sw}=0.0635$. Thm. \ref{thm:main} guarantees that, at least for sufficiently small $\mu$, the solutions are attracted by the target point $\ue{z}^{(1),-}$, during the so-called relaxation stage. In this case, the ``tolerance'' expressed by (\ref{eq:deltaswstar}) corresponds to $\delta_{sw}^*=0.0001367$ which allows an overall error of about $0.2 \%$ around $T_{sw}$ to switch off the field.

\subsection{Case study $2$: ``Experimental''}\label{subsec:two}
\begin{figure*}[t!]\begin{center}
  		\begin{minipage}[c][1\width]{0.49\textwidth}
			\hspace{0pt}
			{\begin{overpic}[width=1\textwidth]{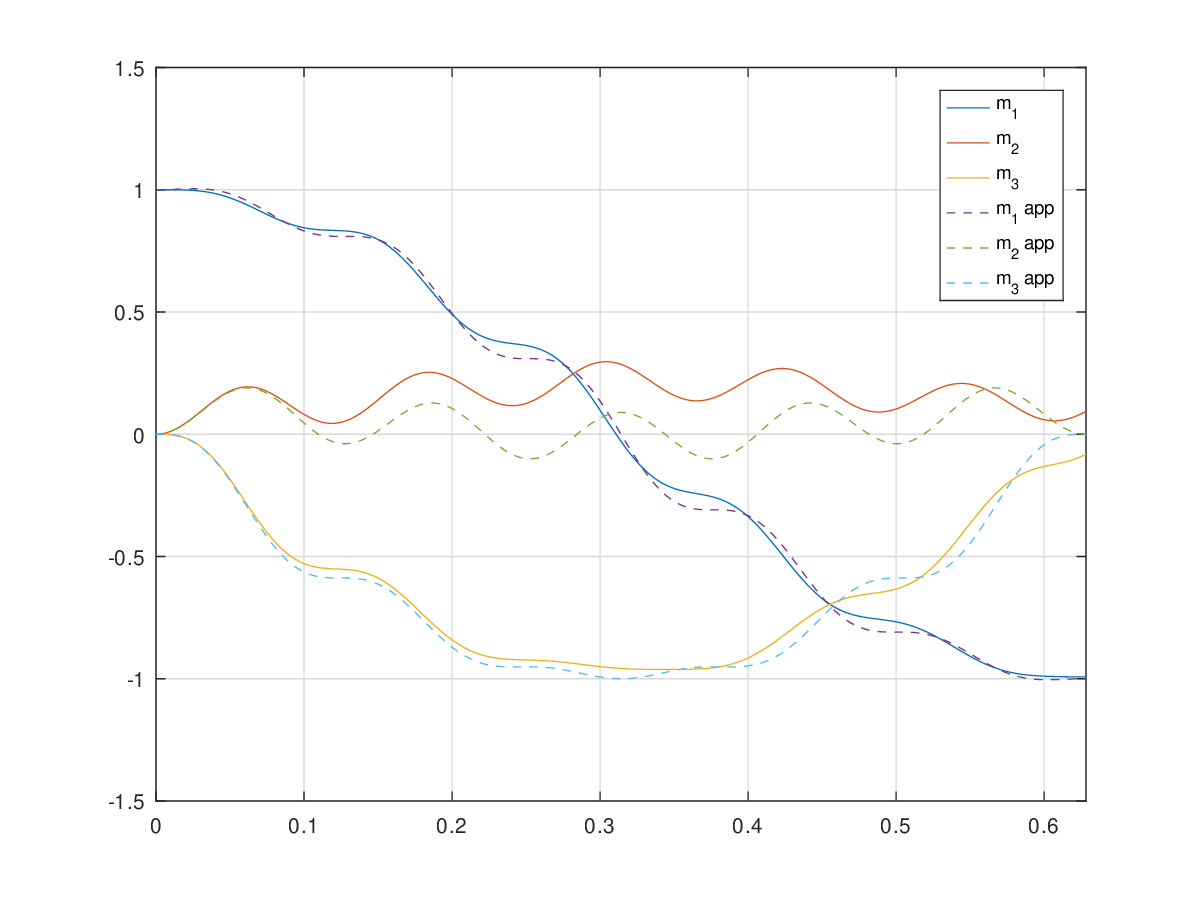}
					\put(0,78){\footnotesize (a)}
			\end{overpic}}
		\end{minipage}	
		\begin{minipage}[c][1\width]{0.49\textwidth}
			\hspace{5pt} 
			{\begin{overpic}[width=\textwidth]{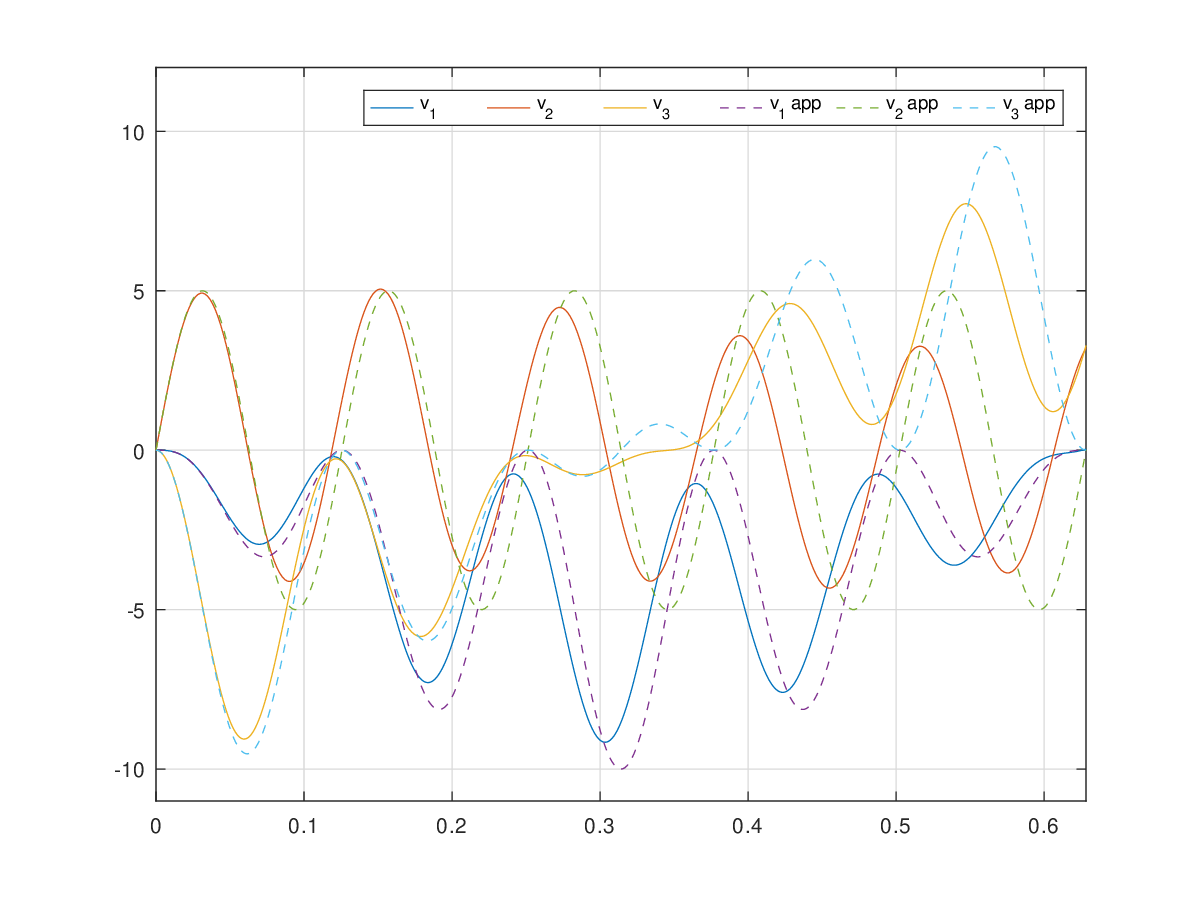}
					\put(0,78){\footnotesize (b)}
			\end{overpic}}
		\end{minipage}\\
		\vspace{-30pt}
  	\begin{minipage}[c][1\width]{0.49\textwidth}
			\hspace{5pt} 
			{\begin{overpic}[width=\textwidth]{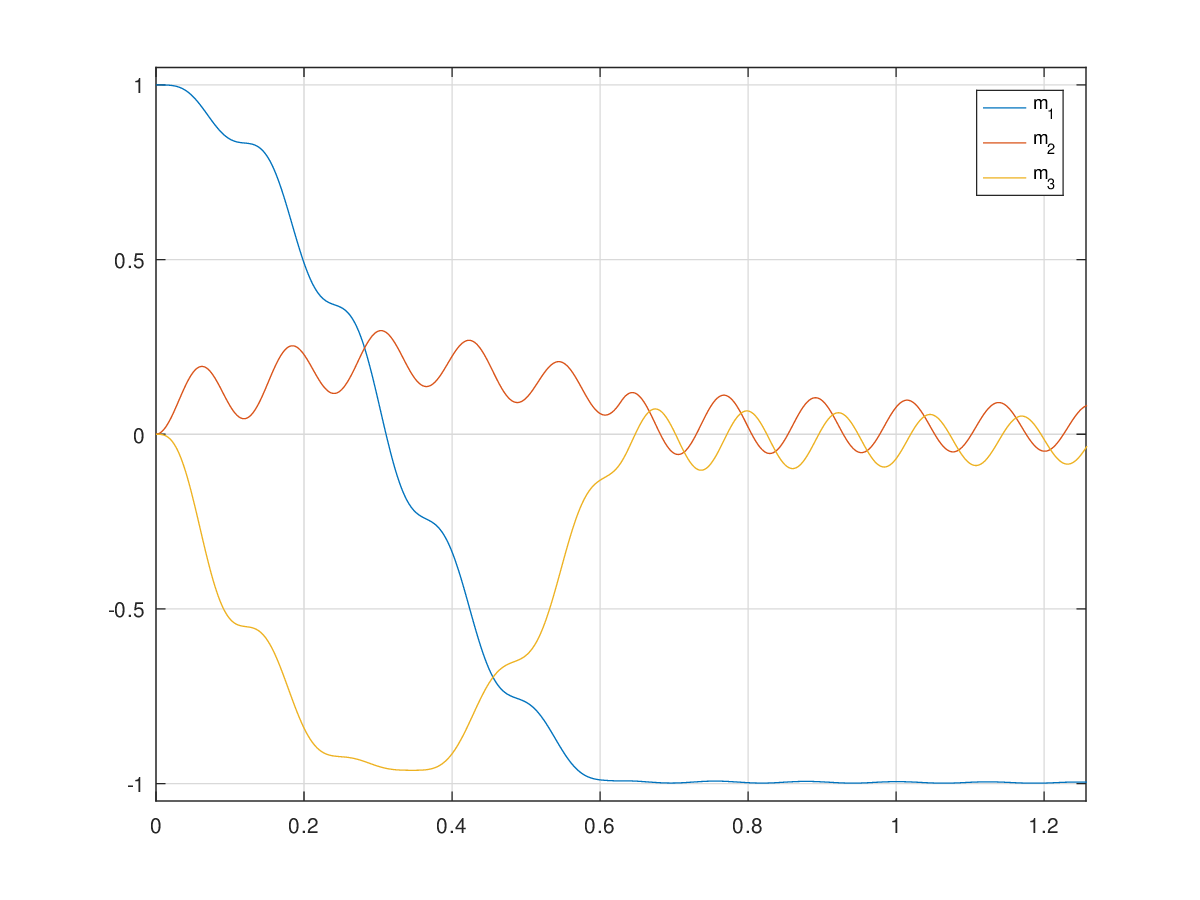}
					\put(0,78){\footnotesize (c)}
			\end{overpic}}
		\end{minipage}	
    \begin{minipage}[c][1\width]{0.49\textwidth}
			\hspace{5pt}
			{\begin{overpic}[width=\textwidth]{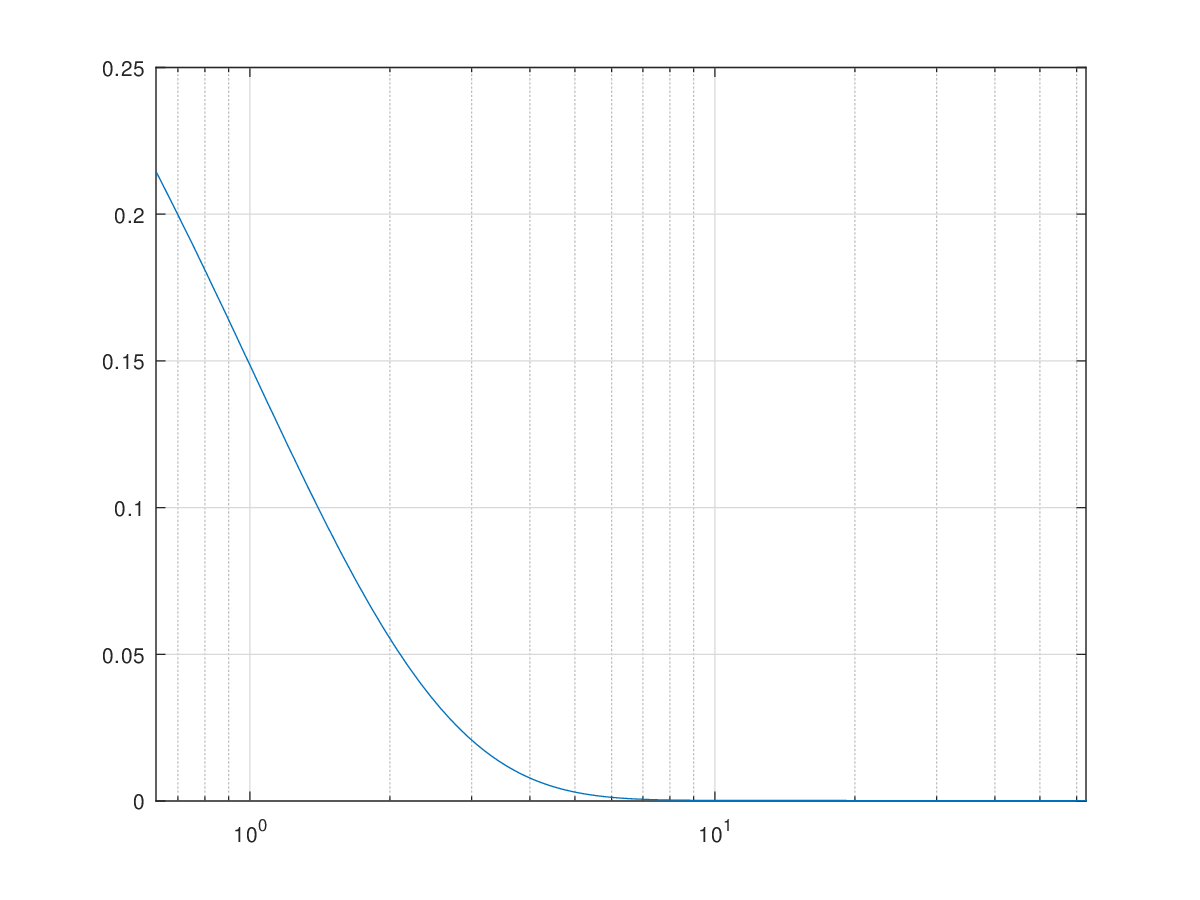}
					\put(0,78){\footnotesize (d)}
 			\end{overpic}}
		\end{minipage}\\
	\end{center}
	\vspace{-30pt}
	\caption{Comparison between the numerically computed solutions (continuous line) of (\ref{eq:illgsixd}) and their approximation provided by (\ref{eq:apporiginalone}) and (\ref{eq:mdot}) (dashed line) for  $\ue{m}(t)$ in panel (a) and  $\ue{v}(t)$ in panel (b). Similarly to Fig. \ref{fig:validation2}, panel (c) shows the behaviour of the computed solution $\ue{m}(t)$ over the time interval $ [0, 2 T_{sw}]$ where $T_{sw}=0.6283$. Finally, the decay of the function $\mathcal{W}$ over the interval $ [T_{sw},100 T_{sw}]$ is reported in  panel (d) in $\log$ scale.\\ Although the approximation is not good as in the examples of sec. \ref{subsec:one}, panel (a) still shows a ``good approximation'' of $m_3(t)$ and, most importantly, $m_1(t)$, i.e. the ``switching'' variable. }
	\label{fig:validation3}
	\vspace{0pt}
\end{figure*}

The aim of this section is to test the validity of Thm. \ref{thm:main}
in a situation relevant to possible experimental observations of inertial switching. In this respect, we refer to the physical parameters used in ref.\cite{neeraj2022inertial}, namely those related to a polycrystalline permalloy extended thin-film with saturation magnetization such that $\mu_0 M_s=0.92$ T and uniaxial anisotropy along the $x$ axis corresponding to an effective field of 100 mT. \\
For this reason we shall start from the values ``without the hats'', i.e. those related to the original model (\ref{eq:illgsixd}). The other quantities will be determined by proceeding ``backwards''.\\
Let us consider the following (paradigmatic) values whose magnitude is fully compatible with experimentally relevant cases
\[
\alpha=0.01 \quad ; \quad \eta=0.02 \mx{,}
\]
while $\uue{D}$ is chosen as in sec. \ref{subsec:one}. \\
To ensure a feasible value for the external field which still provides the required separation between the orders of magnitude of the parameters, we choose $\ep=0.1$. By $(\ref{eq:alphaeta})$ one gets $\hat{\alpha}=1$ and $\hat{\eta}=2$, hence $\mu=0.1$ by (\ref{eq:mu}), i.e. $O(10)$ times bigger than the one chosen in sec. \ref{subsec:one}. For this reason, it is reasonable to expect an overall worsening of the approximation formulae. However, some key features of the approximated solution can be observed also in this case, as it will be shown shortly.\\
By (\ref{eq:defhats}), $\hat{\omega}=2 \omega$. As a general principle, in order to contain the approximation error, $\omega$ should not be too small, as a larger period would require an approximation for longer time intervals. Hence we shall choose $\hat{\omega} \sim 1$, implying that, by (\ref{eq:newmu}), and recalling that $D_{2,1}\sim 0.1$ one finds $\Xi <1$.\\
As a consequence, we need to choose $\hat{\omega}$ as in $ii)$ of Thm. \ref{thm:main}, i.e. $\hat{\omega}=5/n$. Let us choose $n=1$ in such a way $\hat{\omega}=1$ exactly, hence $\omega=0.5$ and $\ue{h}_a=(0,5,0)$. It is easy to check that we still have $\mu \leq \tilde{\mu}_0=1.2299$. \\
A set of numerical experiments similar to those contained in sec. \ref{subsec:one} is reported in Fig. \ref{fig:validation3}.   

One can clearly see (in panel (a)) that the approximated solution (dashed lines) is still in very good agreement with the magnetization component $m_1(t)$ that switches from +1 to -1, which is the most important quantity for describing the switching process (also $m_3(t)$ appears to be close to the approximated solution). The component $m_2(t)$ is the one exhibiting more appreciable deviation from the actual solution, while also the 'velocities' $v_1,v_2,v_3$ (see panel (b) in Fig.\ref{fig:validation3}) seem to be well-captured by the approximated solution for long enough time. It is also apparent (panel (c)) that, by switching the external field off at time $t=T_{sw}$, the magnetization finds itself very close to the target minimum. This means that the system's free energy $\mathcal{W}$ is sufficiently inside the well associated with the reversed magnetization and, therefore, relaxes towards the related minimum (as it is shown in panel (d)) realizing successful switching (i.e. magnetization is 'captured' by the basin of attraction of $m_1=-1$). 

\clearpage

\section*{Appendix A.}

\proof(of Prop. \ref{prop:one})\\
Let us recall the structure of $\Gamma$ in (\ref{eq:Gamma}) and (\ref{eq:sigmaomega}), then define
\beq{eq:m} 
\uue{\mathcal{M}}:= \frac{1}{\sqrt{2} \omega}
\begin{pmatrix}
\sigma & \sigma & \sqrt{2} \hat{h}_1\\
\sigma^{-1} (i \omega \hat{h}_3 - \hat{h}_1 \hat{h}_2) & -\sigma^{-1}( i \omega \hat{h}_3  + \hat{h}_1 \hat{h}_2) & \sqrt{2} \hat{h}_2\\
-\sigma^{-1} (i \omega \hat{h}_2  + \hat{h}_1 \hat{h}_3 )& \sigma^{-1} (i \omega \hat{h}_2 - \hat{h}_1 \hat{h}_3) & \sqrt{2} \hat{h}_3
\end{pmatrix}
\mx{.}
\eeq
The following properties can be proven via a straightforward check: 
$\uue{\mathcal{M}}$ is:
\begin{itemize}
\item unitary, i.e. $\uue{\mathcal{M}}^* \uue{\mathcal{M}} = \uue{I}$, where $\uue{\mathcal{M}}^*$ denotes the conjugate-transpose of $\uue{\mathcal{M}}$ and $\uue{I}$ the unit matrix.
\item diagonalising for $\uue{\Gamma}$, i.e. $\uue{\mathcal{M}}^* \uue{\Gamma} \uue{\mathcal{M}}=\diag(-i \omega, i \omega, 0)$. 
\end{itemize}
Let us now consider a second unitary matrix (straightforward check), defined as 
\beq{eq:u} 
\uue{\mathcal{U}}:= \frac{1}{\sqrt{2}}
\begin{pmatrix}
1 & -i  & 0\\
1 & i  & 0\\
0 & 0 & \sqrt{2} 
\end{pmatrix}
\mx{.}
\eeq
This matrix is a very standard object, in particular, the top-left $2 \times 2$ minor of it arises, for instance, when considering variable transformations from the complex to the real field in the context of Action-Angle variables in Mechanics, see e.g. \cite{gior03}. The proof is complete by setting $\uue{C}:=\uue{\mathcal{M}} \uue{\mathcal{U}}$ and the expression (\ref{eq:c}) is readily obtained. The described properties are easily checked.
\endproof

\section*{Appendix B: Algorithms}

\begin{algorithm}[H]%
\caption{Trajectories approximation}
\label{alg:one}
\begin{algorithmic}
\State \textbf{Input:} $\uue{D}$, $\hat{\ue{h}}_a$ satisfying Hyp \ref{hyp:one}, $\hat{\alpha},\hat{\eta}>0$;
\State $\bullet$ Compute $\sigma,\omega$ via (\ref{eq:sigmaomega});
\State $\bullet$ Construct the matrices $\uue{C}$ via (\ref{eq:c}) then $\uue{E}$ via (\ref{eq:eexplicit});
\State $\bullet$ Determine $\hat{\omega}$ and $\hat{\uue{E}}$ via (\ref{eq:defhats});
\State $\bullet$ Evaluate the threshold $\mu_0$ via (\ref{eq:muzero});
\State $\bullet$ Set an initial condition $\ue{m}(0)$ and compute the corresponding $\theta(0)$ from $\mathrm{S}^{-1} \left(\ue{C}^{-1}\ue{m}(0) \right)$ by using the inverse of (\ref{eq:spherical}); 
\State $\bullet$ Set $\ep>0$ ``sufficiently small'' and compute $\mu$ from (\ref{eq:mu});
\If{$\mu \leq \mu_0$ $\And$ $\theta(0)\neq 0,\pi$} proceed 
\State $\bullet$ Compute the functions $(\theta^{[\leq 2]}(\tau),\ph^{[\leq 2]}(\tau))$ by using (\ref{eq:thphapp}) then use (\ref{eq:spherical}) (with $r\equiv 1$) and finally (\ref{eq:c}) to compute $\ue{m}^{[\leq 2]}(\tau)=\ue{C}  \mathrm{S} (1,\theta^{[\leq 2]}(\tau),\ph^{[\leq 2]}(\tau))$; 
\State $\bullet$ Restore the original time by using (\ref{eq:tau}) to obtain the approximation up to $O(\mu^2)$ of the exact solution of (\ref{eq:illg}) with parameters $\ue{h}_a,\alpha,\eta$ determined via  (\ref{eq:assumptionha}) and (\ref{eq:alphaeta}), respectively;
\State $\bullet$ The solution of (\ref{eq:illgsixd}) is obtained via differentiation w.r.t. $t$ of the already obtained expression. This implies that those will be at least $O(\mu)$ accurate.  
\EndIf  
\end{algorithmic}
\end{algorithm}

\begin{algorithm}[H]%
\caption{Existence of a switching trajectory}
\label{alg:two}
\begin{algorithmic}
\State \textbf{Input:} $\uue{D}$, $\hat{\ue{h}}_a$ satisfying Hyp. \ref{hyp:one} and Hyp. \ref{hyp:one}, $\hat{\alpha},\hat{\eta}>0$;
\State $\bullet$ $\uue{C}$ and $\uue{E}$ are given by 
(\ref{eq:cnew}) and (\ref{eq:newe}), respectively;
\State $\bullet$ Compute $\sigma,\omega,\hat{\omega}$ and $\hat{\uue{E}}$ as in \textbf{Algorithm} \ref{alg:one};  \State $\bullet$ Evaluate the ``new'' threshold $\tilde{\mu}_0$ via (\ref{eq:newmu});
\State $\bullet$ Set $\ep>0$ ``sufficiently small'' and compute $\mu$ from (\ref{eq:mu});
\If{$\mu \leq \tilde{\mu}_0$} proceed 
\State $\bullet$ Compute $\Xi$ via (\ref{eq:newmu});
\State $\bullet$ Compute $K_f$ via (\ref{eq:Kf}) then $\delta_{sw}$ and $T_{sw}$ by using (\ref{eq:tdelta});
\If{$\Xi \geq 1$}
\State $\bullet$ Pick any $t^* \in [T_{sw}-\delta_{sw},T_{sw}+\delta_{sw}]$;  
\State $\bullet$ $\flag \leftarrow \true$;
\ElsIf{$\hat{\omega}$ is of the form $1/(2 \mu n)$ for some $n\in \NN$}
\State $\bullet$ Compute $\delta_{sw}^*$ via (\ref{eq:deltaswstar});
\State $\bullet$ Pick any $t^* \in [T_{sw}-\delta_{sw}^*,T_{sw}+\delta_{sw}^*]$;  
\State $\bullet$ $\flag \leftarrow \true$;
\EndIf
\EndIf
\If{$\flag$}
the system (\ref{eq:illgsixd}) with parameters $\alpha,\eta$ as in (\ref{eq:alphaeta}) and the applied field $\ue{h}_a(t)$ defined by (\ref{eq:switchoff}), undergoes a switching motion in the sense of  Def. \ref{def:switching}.
\EndIf
\end{algorithmic}
\end{algorithm}

\subsection*{Acknowledgements} This work has been supported by the Italian Ministry of University and Research, PRIN2020 funding program, grant number 2020PY8KTC of the University of Naples Federico II. \\
The numerical simulations and the corresponding plots have been performed with GNU Octave \cite{oct}. 

\bibliographystyle{alpha}
\bibliography{cqls.bib}

\newcommand{\etalchar}[1]{$^{#1}$}
\begin{thebibliography}{EBHW20}

\bibitem[BLS{\etalchar{+}}10]{bedau2010spin}
D.~Bedau, H.~Liu, J.Z. Sun, J.A. Katine, E.E. Fullerton, S.~Mangin, and A.D.
  Kent.
\newblock Spin-transfer pulse switching: From the dynamic to the thermally
  activated regime.
\newblock {\em Applied Physics Letters}, 97(26), 2010.

\bibitem[BMSd03]{bertotti2003geometrical}
G.~Bertotti, I.D. Mayergoyz, C.~Serpico, and M.~d'Aquino.
\newblock Geometrical analysis of precessional switching and relaxation in
  uniformly magnetized bodies.
\newblock {\em IEEE transactions on magnetics}, 39(5):2501--2503, 2003.

\bibitem[CRW11]{Ciornei2011magnetization}
M.C. Ciornei, J.M. Rub{\'{\i}}, and J.E. Wegrowe.
\newblock Magnetization dynamics in the inertial regime: Nutation predicted at
  short time scales.
\newblock {\em Physical Review B}, 83(2):020410(R), January 2011.

\bibitem[DPG{\etalchar{+}}20]{Dieny2020}
B.~Dieny, I.~L. Prejbeanu, K.~Garello, P.~Gambardella, P.~Freitas,
  R.~Lehndorff, W.~Raberg, U.~Ebels, S.~O. Demokritov, J.~Akerman, A.~Deac,
  P.~Pirro, C.~Adelmann, A.~Anane, A.~V. Chumak, A.~Hirohata, S.~Mangin,
  Sergio~O. Valenzuela, M.~Cengiz Onbaşlı, M.~d’Aquino, G.~Prenat,
  G.~Finocchio, L.~Lopez-Diaz, R.~Chantrell, O.~Chubykalo-Fesenko, and
  P.~Bortolotti.
\newblock Opportunities and challenges for spintronics in the microelectronics
  industry.
\newblock {\em Nature Electronics}, 3(8):446–459, August 2020.

\bibitem[dPP{\etalchar{+}}23]{daquino2023micromagnetic}
M.~d'Aquino, S.~Perna, M~Pancaldi, R.~Hertel, S.~Bonetti, and C.~Serpico.
\newblock Micromagnetic study of inertial spin waves in ferromagnetic nanodots.
\newblock {\em Physical Review B}, 107(14), April 2023.

\bibitem[dPS24]{DAQUINO2024112874}
M.~d'Aquino, S.~Perna, and C.~Serpico.
\newblock Midpoint geometric integrators for inertial magnetization dynamics.
\newblock {\em Journal of Computational Physics}, 504:112874, 2024.

\bibitem[dSS{\etalchar{+}}04]{daquino2004numerical}
M.~d’Aquino, W.~Scholz, T.~Schrefl, C.~Serpico, and J.~Fidler.
\newblock Numerical and analytical study of fast precessional switching.
\newblock {\em Journal of applied physics}, 95(11):7055--7057, 2004.

\bibitem[EBHW20]{oct}
J.W. Eaton, D.~Bateman, S.~Hauberg, and R.~Wehbring.
\newblock {\em {GNU Octave} version 5.2.0 manual: a high-level interactive
  language for numerical computations}, 2020.

\bibitem[FdS24]{fortunati2024}
Alessandro Fortunati, Massimiliano d’Aquino, and Claudio Serpico.
\newblock Controlled quasi-latitudinal solutions for ultra-fast spin-torque
  magnetization switching.
\newblock {\em International Journal of Bifurcation and Chaos}, 34(05):2450056,
  2024.

\bibitem[FW19]{fliap}
A.~Fortunati and S.~Wiggins.
\newblock {A Lie transform approach to the construction of Lyapunov functions
  in autonomous and non-autonomous systems}.
\newblock {\em Journal of Mathematical Physics}, 60(8):082704, 08 2019.

\bibitem[Gio03]{gior03}
A.~Giorgilli.
\newblock Exponential stability of {H}amiltonian systems.
\newblock In {\em Dynamical systems. {P}art {I}}, Pubbl. Cent. Ric. Mat. Ennio
  Giorgi, pages 87--198. Scuola Norm. Sup., Pisa, 2003.

\bibitem[Har10]{Hartmann2010LectureNO}
C.~Hartmann.
\newblock Lecture notes on singularly perturbed differential equations (edited
  by sebastian john) singularly perturbed differential equations.
\newblock 2010.

\bibitem[Kha02]{khal}
H.K. Khalil.
\newblock {\em Nonlinear Systems}.
\newblock Pearson Education. Prentice Hall, 2002.

\bibitem[KR02]{kaka2002precessional}
S.~Kaka and S.E. Russek.
\newblock Precessional switching of submicrometer spin valves.
\newblock {\em Applied Physics Letters}, 80(16):2958--2960, 2002.

\bibitem[LaS60]{lasalle}
J.~LaSalle.
\newblock Some extensions of {L}iapunov's second method.
\newblock {\em IRE Transactions on Circuit Theory}, 7(4):520--527, 1960.

\bibitem[MBS09]{BMS2009}
I.D. Mayergoyz, G.~Bertotti, and C.~Serpico.
\newblock {\em Nonlinear Magnetization Dynamics in Nanosystems}.
\newblock Elsevier Science, 2009.

\bibitem[NAK{\etalchar{+}}20]{Neeraj2021inertial}
Kumar Neeraj, Nilesh Awari, Sergey Kovalev, Debanjan Polley, Nanna~Zhou
  Hagstr\"{o}m, Sri Sai Phani~Kanth Arekapudi, Anna Semisalova, Kilian Lenz,
  Bertram Green, Jan-Christoph Deinert, Igor Ilyakov, Min Chen, Mohammed
  Bawatna, Valentino Scalera, Massimiliano d'Aquino, Claudio Serpico, Olav
  Hellwig, Jean-Eric Wegrowe, Michael Gensch, and Stefano Bonetti.
\newblock Inertial spin dynamics in ferromagnets.
\newblock {\em Nature Physics}, 17(2):245--250, September 2020.

\bibitem[Nay08]{nay}
A.H. Nayfeh.
\newblock {\em Perturbation Methods}.
\newblock Physics textbook. Wiley, 2008.

\bibitem[NPS{\etalchar{+}}22]{neeraj2022inertial}
K.~Neeraj, M.~Pancaldi, V.~Scalera, S.~Perna, M.~d'Aquino, C.~Serpico, and
  S.~Bonetti.
\newblock Magnetization switching in the inertial regime.
\newblock {\em Physical Review B}, 105:054415, Feb 2022.

\bibitem[SdBM09]{serpico2009analytical}
C.~Serpico, M.~d'Aquino, G.~Bertotti, and I.D. Mayergoyz.
\newblock Analytical description of quasi-random magnetization relaxation to
  equilibrium.
\newblock {\em IEEE transactions on magnetics}, 45(11):5224--5227, 2009.

\end{thebibliography}

\end{document}